\documentclass[11pt,twoside]{amsart}

\usepackage[dvipsnames]{xcolor}
\usepackage{hyperref}
\hypersetup{colorlinks=true, linkcolor=Blue, citecolor=Maroon}
\usepackage{amsthm, amsmath, amscd, amssymb,centernot,txfonts}
\usepackage{tikz-cd}
\usepackage{float}
\usepackage{multirow}
\usepackage{lipsum}
\usepackage{multicol}
\usepackage{amsmath}
\usepackage{amssymb}

\usepackage[normalem]{ulem}
\usepackage{amsfonts}
\usepackage{mathrsfs}
\usepackage{parskip}
\usepackage[all]{xy}
\usepackage[left=2.5cm,top=2.5cm,bottom=3cm,right=2.5cm]{geometry}
\usepackage{dirtytalk}
\usepackage{verbatim}
\usepackage{mathtools}
\usetikzlibrary{shapes,arrows}
\usepackage{etoolbox}
\usepackage{dynkin-diagrams}

\usepackage{enumitem}
\usepackage{chngpage}

\usepackage{adjustbox}
\usepackage[utf8]{inputenc}
\usepackage{fourier} 
\usepackage{array}
\usepackage{makecell}

\usepackage[linktocpage=true]{}

\usepackage{color, colortbl}
\definecolor{LightCyan}{rgb}{0.88,1,1}

\usepackage[normalem]{ulem}

\usepackage[first=0,last=9]{lcg}

\newcommand{\stkout}[1]{\ifmmode\text{\sout{\ensuremath{#1}}}\else\sout{#1}\fi}

\theoremstyle{plain}
\numberwithin{equation}{section}
\newtheorem{theorem}{Theorem}[section]
\newtheorem{proposition}[theorem]{Proposition}
\newtheorem{lemma}[theorem]{Lemma}
\newtheorem{corollary}[theorem]{Corollary}

\newtheorem{question}[theorem]{Question}

\newtheorem{set-up}[theorem]{Set-up}

\theoremstyle{definition}
\newtheorem{remark}[theorem]{Remark}

\newtheorem{definition}[theorem]{Definition}

\newcommand*{\QEDB}{\hfill\ensuremath{\square}}

\tikzstyle{decision} = [diamond, draw, , 
    text width=4.5em, text badly centered, node distance=3cm, inner sep=0pt]
\tikzstyle{block} = [rectangle, draw, , 
    text width=10em, text centered, rounded corners, minimum height=2em]
\tikzstyle{block1} = [rectangle, draw, , 
    text width=5em, text centered, rounded corners, minimum height=2em]    
\tikzstyle{line} = [draw, -latex']
\tikzstyle{cloud} = [draw, ellipse,, node distance=3cm,
    minimum height=2em]

\begin{document}
\title [Koszul property of Ulrich bundles]{Koszul property of Ulrich bundles and rationality of moduli spaces of stable bundles on Del Pezzo surfaces}

\author[P. Bangere]{Purnaprajna Bangere}
\address{Department of Mathematics, University of Kansas, Lawrence, USA}
\email{purna@ku.edu}

\author[J. Mukherjee]{Jayan Mukherjee}
\address{Department of Mathematics, Oklahoma State University, Stillwater, USA}
\email{mukherjeejayan@gmail.com}

\author[D. Raychaudhury ]{Debaditya Raychaudhury} \thanks{Corresponding author: Debaditya Raychaudhury, email id: rcdeba@gmail.com}
\address{Department of Mathematics, University of Toronto, Toronto, Canada}
\email{debaditya.raychaudhury@utoronto.ca}
\address{Current address of D. Raychaudhury: Department of Mathematics, University of Arizona, Tucson, USA}
\email{draychaudhury@math.arizona.edu, rcdeba@gmail.com}


\maketitle

\begin{abstract}
Let $\mathcal{E}$ be a vector bundle on a smooth projective variety $X\subseteq\mathbb{P}^N$ that is Ulrich with respect to the hyperplane section $H$. In this article, we study the Koszul property of $\mathcal{E}$, the slope--semistability of the $k$--th iterated syzygy bundle $\mathcal{S}_k(\mathcal{E})$ for all $k\geq 0$ and rationality of moduli spaces of slope--stable bundles on Del Pezzo surfaces. As a consequence of our study, we show that if $X$ is a Del Pezzo surface of degree $d\geq 4$, then any Ulrich bundle $\mathcal{E}$ satisfies the Koszul property and is slope--semistable. We also show that, for infinitely many Chern characters ${\bf v}=(r,c_1, c_2)$, the corresponding moduli spaces of slope--stable bundles $\mathfrak{M}_H({\bf v})$ when non--empty, are rational, and thereby produce new evidences for a conjecture of Costa and Mir\'o-Roig. As a consequence, we show that the iterated syzygy bundles of Ulrich bundles are dense in these moduli spaces. 
\end{abstract}
\section{Introduction} 

Let $X\subseteq\mathbb{P}^N$ be a smooth projective variety of dimension $n$ and let $\mathcal{E}$ be an arithmetically Cohen--Macaulay vector bundle on $X$. Set $\Gamma_*(\mathcal{E}):=\bigoplus_{m\geq 0}H^0(\mathcal{E}(m))$ and notice that $\Gamma_*(\mathcal{E})$ is a graded Cohen-Macaulay module over $S:=\textrm{Sym}(H^0(\mathcal{O}_{\mathbb{P}^N}(1)))$. The vector bundle $\mathcal{E}$ is called an {\it Ulrich bundle} if $\Gamma_*(\mathcal{E})$ admits a linear minimal graded free resolution as an $S$ module i.e., a minimal resolution of the following form:
$$0\to S(-(N-n))^{\oplus \beta_{N-n}}\to S(-(N-n-1))^{\oplus \beta_{N-n-1}}\to\cdots\to S^{\beta_0}\to \Gamma_*(\mathcal{E})\to 0.$$ 
Graded Cohen--Macaulay modules $M$ that admit a minimal graded free resolution as above are called {\it Ulrich modules} (also known as {\it linear maximal Cohen--Macaulay modules}); they were first introduced in commutative algebra by Ulrich in \cite{Ulr}. In their seminal works, Beauville in \cite{BeaU2}, and Eisenbud, Schreyer and Weyman in \cite{ES} studied them in the context of algebraic geometry. It follows from \cite{ES} that a vector bundle $\mathcal{E}$ on $X$ is Ulrich if and only if $H^{\bullet}(\mathcal{E}(-p))=0$ for all $1\leq p\leq \dim X$. 

\smallskip

\noindent{\bf Background.} The study of Ulrich bundles is an active area of research in the present years. For example the existence questions for these bundles are studied by many authors; see  \cite{ACC}, \cite{ACM}, \cite{AFO}, \cite{BeaU1}, \cite{Cas3}, \cite{Cas1}, \cite{Cas2},\cite{CH2}, \cite{Fae}, \cite{Jos}, \cite{Lop2}, \cite{Lop3}, \cite{PT}, \cite{RT} and the references therein. The relations between these bundles with determinantal representations, Chow forms, Clifford algebras were studied in \cite{BeaU2}, \cite{ES}, \cite{CKM1}, \cite{CKM2}. Ulrich bundles on cubic surfaces were studied by Casanellas and Hartshorne in \cite{CHGS}, and the problem of determining the determinant of Ulrich bundles on arbitrary Del Pezzo surfaces has been connected to the Minimal Resolution Conjecture by the work of Coskun, Kulkarni and Mustopa in \cite{CKM} (see also \cite{MiP}, \cite{MiP1}). The positivity of these bundles have also been investigated by Lopez, Sierra and Mu\~{n}oz in \cite{LC}, \cite{LM}, \cite{Lop1}. We refer to \cite{BeaU} for a beautiful introduction to  Ulrich bundles.

\smallskip

\noindent{\bf Objective.} The main aim for this article is to study the syzygy bundles of Ulrich bundles. We recall that for any globally generated vector bundle $\mathcal{E}$ on a smooth projective variety $X$, the {\it syzygy bundle} (also known in literature as the {\it kernel bundle}, the {\it dual span bundle}, or the {\it Lazarsfeld--Mukai bundle}) $M_{\mathcal{E}}$ is defined as the kernel of the evaluation map $H^0(\mathcal{E})\otimes\mathcal{O}_X\to \mathcal{E}$, i.e., it fits into the following exact sequence 
\begin{equation}\label{defining}
    0\to M_{\mathcal{E}}\to H^0(\mathcal{E})\otimes\mathcal{O}_X\to \mathcal{E}\to 0.
\end{equation}

It was shown by Casanellas and Hartshorne that Ulrich vector bundles on $X\subseteq\mathbb{P}^N$ are Gieseker--semistable with respect to the hyperplane bundle $H$, and further the notions of Gieseker and slope--stability with respect to $H$ coincides for them (see \cite[Theorem 2.9]{CHGS}). It is thus natural to study the semi--stability for the syzygy bundles associated to these bundles and we give a sufficient condition for slope--semistability (with respect to $H$) based on the intersection numbers $H^n$ and $H^{n-1}\cdot K_X$ in Proposition \ref{semistability}. The criterion described in this proposition is satisfied for any embedding of anticanonical rational surfaces (see \cite[Theorem III.1]{Har}), in particular for anticanonically embedded Del Pezzo surfaces (by which we mean blow up of the projective plane in $9-d$ points in general position).

\smallskip

Even though the syzygy bundle $M_{\mathcal{E}}$ itself has no global sections, one can still try to find out when $M_{\mathcal{E}}(H)$ is globally generated. Further, if $M_{\mathcal{E}}(H)$ is globally generated, then we can analogously study the syzygy bundle of $\mathcal{S}_0(\mathcal{E}):=M_{\mathcal{E}}(H)$. We define a property $(\textrm{Kos}_k)$ for a vector bundle $\mathcal{E}$ on $X$ (dependent on $H$) in Subsection \ref{33}. It follows from our definition that if a vector bundle $\mathcal{E}$ satisfies this property with respect to $H$, then we can iterate the construction and define $\mathcal{S}_p(\mathcal{E})$ for all $0\leq p\leq k-1$. In particular, if $\mathcal{E}$ satisfies property $(\textrm{Kos}_{\infty})$ then $\mathcal{S}_{p}(\mathcal{E})$ exists for all $p\geq 0$. Recall that for a graded module $M$ over a graded $\mathbb{C}$--algebra $R:=\oplus_{i\geq 0} R_i$ with $R_0=\mathbb{C}$ and $R$ generated in degree one, we say that $M$ is {\it Koszul} if $M$ admits a linear minimal graded free resolution i.e., a minimal resolution of the form 
$$\cdots\to E_{i+1}\to E_i\to\cdots\to E_1\to E_0\to M\to 0$$ where $E_i=R(-i)^{\beta_i}$. 
An argument of Lazarsfeld (which follows an interpretation of Kempf) as described by Pareschi in \cite{Par}, shows that if $\mathcal{E}$ is Ulrich on an arithmetically Cohen--Macaulay variety $X\subseteq\mathbb{P}^N$ and if $\mathcal{S}_p(\mathcal{E})$ exists for all $p\geq 0$, then $\mathcal{E}$ satisfies property $(\textrm{Kos}_{\infty})$ if and only if $\Gamma_*(\mathcal{E})$ is a Koszul $R(H)$ module where $R(H):=S/I$ is the coordinate ring of $X\subseteq\mathbb{P}^N$, (i.e., $S:=\textrm{Sym}(H^0(\mathcal{O}_{\mathbb{P}^N}(1)))$ and $I:=I_{X/\mathbb{P}^N}$ is the homogeneous ideal of the embedding) when $X$ is embedded as a projectively normal variety (see Corollary \ref{kosiff2}). Further, it turns out (see Theorem \ref{kosiff}) that, for an Ulrich vector bundle $\mathcal{E}$ (with respect to the hyperplane bundle) on an arithmetically Cohen--Macaulay and projectively normal variety $X\subseteq\mathbb{P}^N$, $\Gamma_*(\mathcal{E})$ is a Koszul $R(H)$ module if and only if the embedding $X\subseteq\mathbb{P}^N$ is Koszul, i.e., $R(H)$ is a {\it Koszul $\mathbb{C}$--algebra} (which means that $\mathbb{C}$ as a graded module concentrated in degree zero is Koszul when viewed as a module over the homogeneous coordinate ring $R(H)$ of the embedding). Now, for an anticanonically embedded Del Pezzo surface of degree $d\geq 4$, we show that any Ulrich bundle $\mathcal{E}$ satisfies $(\textrm{Kos}_{\infty})$ and the iterated syzygy bundles $\mathcal{S}_p(\mathcal{E})$ are 
slope--semistable  with respect to the hyperplane (which is same as anticanonical) bundle (see Theorem \ref{mains}, Corollary \ref{extis}). 

\textcolor{black}{We show that for} an anti--canonically embedded Del Pezzo surfaces $X_d$ with $d\geq 4$, 
the set--theoretic association $\mathcal{F}\mapsto\mathcal{S}_p(\mathcal{F})$ defines a birational map between the moduli spaces of slope--stable bundles (with respect to the hyperplane $=$ anticanonical bundle) with Chern characters ${\bf v}(\mathcal{E})$ and ${\bf v}(\mathcal{S}_p(\mathcal{E}))$ where $\mathcal{E}$ is a stable Ulrich bundle of rank $\geq 2$ whenever the latter space is non--empty. 
Combining this result with the work of Costa and Mir\'o-Roig, we prove the rationality of infinitely many moduli spaces of slope--stable bundles for Del Pezzo surfaces. For $d=3$, we prove analogous birationality result only when $p=0$. \textcolor{black}{In this case, we note that the bundle $\mathcal{E}$ is not Koszul (see Remark \ref{rmkc}).}

\smallskip

\noindent{\bf Main results.} We prove the following theorem in this article (see Subsection \ref{42} for notation).

\begin{theorem}\label{main}
Let $X_d\subseteq\mathbb{P}^d$ be the anti--canonically embedded Del Pezzo surface of degree $3\leq d\leq 8$. Let $\mathcal{E}$ be a stable Ulrich vector bundle with respect to $H$ of rank $r$ where $H:=\mathcal{O}_{\mathbb{P}^d}(1)|_{X_d}=-K_{X_d}$. 
\begin{itemize}
    \item[(1)] If $d=3$ then $M_{\mathcal{E}}$ is simple and slope--semistable with respect to $H$. Moreover, the moduli spaces of slope--stable bundles  $\mathfrak{M}_H({\bf v}(\mathcal{E}))$ and $\mathfrak{M}_H({\bf v}(M_{\mathcal{E}}))$ are both non--empty, smooth, irreducible and are birational to each other. Also, the general bundle in the latter moduli space is of the form $M_{\mathcal{F}}$ for some $\mathcal{F}\in \mathfrak{M}_{H}({\bf v}(\mathcal{E}))$.
    \item[(2)] If $d\geq 4$ then $\mathcal{S}_p(\mathcal{E})$ is simple and slope--semistable with respect to $H$ for all $p\geq 0$. Moreover, for any $p\geq 0$, the moduli space of slope--stable bundles $\mathfrak{M}_H({\bf v}(\mathcal{E}))$ and $\mathfrak{M}_H({\bf v}(\mathcal{S}_{p}(\mathcal{E})))$ are both smooth, irreducible and are birational to each other whenever the latter space is non--empty; in this case the general bundle in the latter moduli space is of the form $\mathcal{S}_p(\mathcal{F})$ for some $\mathcal{F}\in \mathfrak{M}_{H}({\bf v}(\mathcal{E}))$.
\end{itemize}
\end{theorem}

For a smooth projective polarized surface, a moduli space $\mathfrak{M}_H({\bf v})$ of slope--stable vector bundles with respect to $H$ is said to satisfy the {\it weak Brill--Noether property} if a general bundle in the moduli space has no cohomology. This property is useful to show that the moduli space $\mathfrak{M}_H^{ss}({\bf v})$ of Gieseker semistable sheaves admits an effective theta divisor $$\Theta:=\left\{\mathcal{F}\in \mathfrak{M}^{ss}_H({\bf v}) \,\,|\,\, H^1(\mathcal{F})\neq 0\right\}.$$ Weak Brill--Noether property for moduli spaces on rational surfaces were studied by Coskun and Huizenga in \cite{CH1}, \cite{CH2}). As an immediate consequence of our theorem, we obtain the following

\begin{corollary}\label{main3}
Let $X_d\subseteq\mathbb{P}^d$ be an anti--canonically embedded Del Pezzo surface of degree $d$. We set $H:=\mathcal{O}_{\mathbb{P}^d}(1)|_{X_d}=-K_{X_d}$. Let $\mathcal{E}$ be a stable Ulrich vector bundle with respect to $H$ of rank $r$. 
\begin{itemize}
    \item[(1)] When $d=3$, the moduli space $\mathfrak{M}_{H}({\bf v}(M_{\mathcal{E}}))$ of slope--stable vector bundles with respect to $H$ satisfies weak Brill--Noether property.
    \item[(2)] Same conclusion holds for $\mathfrak{M}_{H}({\bf v}(\mathcal{S}_{p}(\mathcal{E})(-H)))$ when $4\leq d\leq 8$, whenever the space is non--empty.
\end{itemize}
\end{corollary}

Now we discuss a third consequence of Theorem \ref{main} concerning rationality. The Unirationality of the moduli spaces of Gieseker--semistable sheaves were studied by Coskun and Huizenga in \cite{CH2}. Further, the rationality of the moduli spaces of slope--stable rank two bundles has been studied in a series of papers by Costa and Mir\'o-Roig (see \cite{CM1}, \cite{CM3}, \cite{CM2}). In particular, they completely solved the problem for Del Pezzo surfaces and established strong partial results for higher rank bundles. Costa and Mir\'o-Roig also conjectured that the moduli spaces of slope--stable bundles on smooth rational surfaces are rational (see \cite{CM1}). Combining our result with their result for rank two bundles, we show rationality for infinitely many moduli spaces of higher rank slope--stable vector bundles on Del Pezzo surfaces $X_d$ for $4\leq d\leq 7$, all of which have dimensions either $1,3$, or $5$. More precisely, for $d=5,6,7$, set $$\alpha_{d,1}:=(1/2)((d-2)+\sqrt{d(d-4)})\quad \textrm{and}\quad \alpha_{d,2}:=(1/2)((d-2)-\sqrt{d(d-4)}).$$ 
For $k\geq -1$ and $4\leq d\leq 7$, let $N_{d,k}$ be the sequence of integers and $c_{i,\{d,k\}}$ be the Chern classes on anticanonically embedded Del Pezzo surfaces $X_d\subseteq\mathbb{P}^d$ with $4\leq d\leq 7$ as described below:
\begin{equation}\label{c0f}
N_{d,k} = 
\begin{cases}
    4k+6, & \text{for } d=4, \\
    \frac{2\left(\left(\alpha_{d,2}^{-(k+2)}+\alpha_{d,2}^{-(k+1)}\right)-\left(\alpha_{d,1}^{-(k+2)}+\alpha_{d,1}^{-(k+1)}\right)\right)}{\sqrt{d(d-4)}} & \text{for } d=5,6,7; \\
\end{cases}
\end{equation}

\begin{equation}\label{c1f}
    \begin{gathered}
    c_{1,\{d,-1\}}={\bf c}_{1,\{d,-1\}},\quad c_{1,\{d,0\}}=-{\bf c}_{1,\{d,-1\}},\\
    c_{1,\{d,k\}}=(-1)^{k+1}{\bf c}_{1,\{d,-1\}}+\sum_{i=0}^{k-1}(-1)^{k+i}N_{d,i}H_d\, \textrm{ for all $k\geq 1$};
\end{gathered}
\end{equation}

\begin{equation}\label{c2f}
\begin{gathered}
c_{2,\{d,-1\}}={\bf c}_{2,\{d,-1\}},\quad c_{2,\{d,0\}}=({\bf c}_{1,\{d,-1\}})^2-{\bf c}_{2,\{d,-1\}},\\
    \begin{array}{ll}
        c_{2,\{d,k\}}=  & \sum\limits_{i=0}^{k-1}(-1)^{k+i+1}\left(1-\binom{N_{d,i}}{2}\right)(c_{1,\{d,i\}})^2+\sum\limits_{i=0}^{k-1} (-1)^{k+i+1}(N_{d,i}+1)c_{1,\{d,i\}}H_d+\\
        & \sum\limits_{i=0}^{k-1}dN_{d,i}^2+(-1)^k({\bf c}_{1,\{d,-1\}}^2-{\bf c}_{2,\{d,-1\}})\, \textrm{ for all $k\geq 1$}.
     \end{array}
\end{gathered} 
\end{equation}

\begin{theorem}\label{main2}
Let $\mathcal{E}$ be a stable Ulrich vector bundle (with respect to $H=-K_{X_d}$) of rank two on $X_d$. Set $${\bf v}_{d,k}:=(N_{d,k},c_{1,\{d,k\}}, c_{2,\{d,k\}}),\quad{\bf c}_{i,\{d,-1\}}=c_i(\mathcal{E})\textrm{ for $i=1,2$.}$$
where $N_{d,k}, c_{1,\{d,k\}}, c_{2,\{d,k\}}$ are as described in \eqref{c0f}, \eqref{c1f} and \eqref{c2f} respectively.
If $c_1(\mathcal{E})^2$ (and $c_2(\mathcal{E})$) takes the values described in the following table, then the moduli space of slope--stable vector bundles $\mathfrak{M}_{H}({\bf v}_{d,k})$ is non--empty smooth irreducible and rational for $k\gg 0$ with dimensions as indicated below.
\begin{table}[H]
    \centering
    \begin{tabular}{  c | c| c | c   }
 \hline
$d$ & ${\bf c}_{1,\{d,-1\}}^2$ & ${\bf c}_{2,\{d,-1\}}$ & $\dim \mathfrak{M}_{H}({\bf v}_{d,k})$ \\
\hline\hline
\multirow{2}{*}{\centering $d=4$} & $12$ & $4$ & $1$\\
\cline{2-4} & $16$ & $6$ & $5$\\
\hline\hline
\multirow{2}{*}{\centering $d=5$} & $16$ & $5$ & $1$\\
\cline{2-4} & $20$ & $7$ & $5$\\
\hline\hline
\multirow{2}{*}{\centering $d=6$} & $20$ & $6$ & $1$\\
\cline{2-4} & $24$ & $8$ & $5$\\
\hline\hline
\multirow{3}{*}{\centering $d=7$} & $24$ & $7$ & $1$\\
\cline{2-4} & $26$ & $8$ & $3$\\
\cline{2-4} & $28$ & $9$ & $5$\\
\hline
\end{tabular}
\end{table}
\end{theorem}

We will also prove rationality of certain moduli spaces of slope--stable bundles for cubic surfaces whose points arise as syzygy bundles 
of stable Ulrich bundles (see Corollary \ref{corcubic}). Two remarks are in order. The first one is that the rank two stable Ulrich bundles on anticanonically embedded Del Pezzo surfaces were classified by Casnati in \cite{Cas2}. In particular, it follows from \cite[Proposition 5.1, Examples 6.4 and 6.5]{Cas2} that there exists stable rank two Ulrich bundles with the values $c_1^2$ described in the theorem above. Second remark is that the smoothness and the irreducibility part of the above result is a theorem of Walter (see \cite{Wal}), and the non--emptiness part is a fundamental theorem on moduli of vector bundles on surfaces (see \cite{O'Gr}). 

\smallskip

\noindent{\bf Organization.} Now we give the structure of this article. In Section \ref{2}, we compare the dimensions of $\textrm{Ext}^i(M_{\mathcal{E}},M_{\mathcal{E}})$ and $\textrm{Ext}^i(\mathcal{E},\mathcal{E})$ for a globally generated bundle $\mathcal{E}$ that we will require in order to prove the main theorem on the birationality of the moduli spaces. In Section \ref{3}, we first study the slope--semistability of the syzygy bundles of Ulrich bundles. We then study the Koszul property for these bundles and study the slope--semistability for the iterated syzygy bundles. We prove the main results stated in the introduction in Section \ref{4}. 

\smallskip

\noindent{\bf Notation and conventions.} We will work over the field of complex numbers $\mathbb{C}$. A variety is an integral separated scheme of finite type over $\mathbb{C}$. For a smooth projective variety $X$, $K_X$ denotes its canonical bundle. For a vector bundle $\mathcal{F}$ on a smooth projective surface $X$, by ${\bf v}(\mathcal{F})$, we denote the triple ${\bf v}(\mathcal{F}):=(\textrm{rank}(\mathcal{F}), c_1(\mathcal{F}), c_2(\mathcal{F}))$.

\smallskip

\noindent{\bf Acknowledgements.} We are grateful to Professor Laura Costa for pointing out an error in an earlier draft of this article, and to Professor Yusuf Mustopa for suggesting several improvements. 
We thank Professor Srikanth Iyengar for helpful and motivating discussions. We are also thankful to the referee for several valuable comments and for suggesting us the use of the diagram \eqref{referee}. During this work,  
the second author was supported by the National Science Foundation, Grant No. DMS-1929284 while in residence at the Institute for Computational and Experimental Research in Mathematics in Providence, RI, as part of the ICERM Bridge program. The research of the third author was supported by a Simons Postdoctoral Fellowship from the Fields Institute for Research in Mathematical Sciences.

\section{The relations between \texorpdfstring{$\textrm{Ext}^i(M_{\mathcal{E}},M_{\mathcal{E}})$}{TEXT} and \texorpdfstring{$\textrm{Ext}^i(\mathcal{E},\mathcal{E})$}{TEXT}}\label{2}

Throughout this section, unless otherwise stated, $X$ is a smooth projective variety of dimension $n\geq 2$, $\mathcal{E}$ is a globally generated vector bundle on $X$ of rank$(\mathcal{E})=r$. Observe that it follows from \eqref{defining} that 
\begin{equation}\label{cime}
    c_1(M_{\mathcal{E}})=-c_1(\mathcal{E}),\quad c_2(M_{\mathcal{E}})=(c_1(\mathcal{E}))^2-c_2(\mathcal{E}).
\end{equation}

\begin{lemma}\label{elt}
$H^0(M_{\mathcal{E}})=0$. Further, the following statements hold
\begin{itemize}
    \item[(1)] if $H^1(\mathcal{O}_X)=0$ then $H^1(M_{\mathcal{E}})=0$,
    \item[(2)] if $\mathcal{O}_X$ has no higher cohomology, then $H^{i+1}(M_{\mathcal{E}})\cong H^i(\mathcal{E})$ for all $i\geq 1$. 
\end{itemize}
\end{lemma}

\noindent\textit{Proof.} Immediately follows from the long exact sequence of cohomology associated to \eqref{defining}.\QEDB

\vspace{5pt}

The goal of this section is to relate the spaces $\textrm{Ext}^i(M_{\mathcal{E}},M_{\mathcal{E}})\cong H^i(M_{\mathcal{E}}\otimes M_{\mathcal{E}}^*)$ and $\textrm{Ext}^i(\mathcal{E},\mathcal{E})\cong H^i(\mathcal{E}\otimes\mathcal{E}^*)$ to compare the deformation theory of $M_{\mathcal{E}}$ against that of $\mathcal{E}$.

\begin{proposition}\label{hom}
Assume $H^{n-1}(K_X\otimes\mathcal{E})=0$ and $H^{n}(K_X\otimes\mathcal{E})=0$. Then,
\begin{itemize}
    \item[(1)] $H^0(M_{\mathcal{E}}^*)\cong H^0(\mathcal{E})^*$,
    \item[(2)] Hom$(M_{\mathcal{E}},M_{\mathcal{E}})\cong\textrm{Hom}(\mathcal{E},\mathcal{E})$,
    \item[(3)] $\mathcal{E}\cong \left(M_{M_{\mathcal{E}}^*}\right)^*$.
\end{itemize}
\end{proposition}

\noindent\textit{Proof.} The arguments below are inspired by the proof of \cite[Proposition 1]{Cam}.

(1) Dualizing the defining exact sequence \eqref{defining} of $M_{\mathcal{E}}^*$, we obtain the following exact sequence
\begin{equation}\label{dualdef}
    0\to \mathcal{E}^*\to H^0(\mathcal{E})^*\otimes \mathcal{O}_X\to M_{\mathcal{E}}^*\to 0.
\end{equation}
By our hypothesis, using Serre duality we obtain $H^0(\mathcal{E}^*)=H^1(\mathcal{E}^*)=0$. Consequently, the long exact sequence of cohomology associated to \eqref{dualdef} proves the claim.

(2) The sequences \eqref{defining}, \eqref{dualdef} yield the following commutative diagram with exact rows and columns:
\begin{equation}\label{referee}
\begin{tikzcd}
    & 0\arrow[d] & 0\arrow[d] & 0\arrow[d] &\\
    0 \arrow{r} & \mathcal{E}^*\otimes M_{\mathcal{E}} \arrow{r} \arrow[d] & H^0(\mathcal{E})\otimes\mathcal{E}^* \arrow{r} \arrow{d} & \mathcal{E}\otimes \mathcal{E}^*\arrow{r}\arrow[d] & 0\\
   0 \arrow{r} & H^0(\mathcal{E})^*\otimes M_{\mathcal{E}} \arrow{r}\arrow[d] & H^0(\mathcal{E})\otimes H^0(\mathcal{E})^*\otimes\mathcal{O}_{X} \arrow{r}\arrow[d] & \mathcal{E}\otimes H^0(\mathcal{E})^*\arrow{r}\arrow[d, "\eta"] & 0\\
   0\arrow[r] & M_{\mathcal{E}}\otimes M_{\mathcal{E}}^*\arrow[r]\arrow[d] & H^0(\mathcal{E})\otimes M_{\mathcal{E}}^*\arrow[r, "\xi"]\arrow[d] & \mathcal{E}\otimes M_{\mathcal{E}}^*\arrow[r]\arrow[d] & 0\\
   & 0 & 0 & 0 &
\end{tikzcd} 
\end{equation}
Consequently, passing to the cohomology of the above, we obtain the following diagram whose left vertical map is an isomorphism since $H^0(\mathcal{E}^*)=H^1(\mathcal{E}^*)=0$:
\begin{equation}\label{cohreferee}
    \begin{tikzcd}
     H^0(\mathcal{E})\otimes H^0(\mathcal{E})^* \arrow[r, "\cong"]\arrow[d, "\cong"] & H^0(\mathcal{E})\otimes H^0(\mathcal{E})^*\arrow[d, "H^0(\eta)"]\\
   H^0(\mathcal{E})\otimes H^0(M_{\mathcal{E}}^*)\arrow[r, "H^0(\xi)"] & H^0(\mathcal{E}\otimes M_{\mathcal{E}}^*)
\end{tikzcd}
\end{equation}
Observe that we also have the isomorphisms $\textrm{Ker}(H^0(\eta))\cong H^0(\mathcal{E}\otimes\mathcal{E}^*)$ and $\textrm{Ker}(H^0(\xi))\cong H^0(M_{\mathcal{E}}\otimes M_{\mathcal{E}}^*)$ which prove the assertion.

(3) Observe that by \eqref{dualdef}, $M_{\mathcal{E}}^*$ is globally generated. 
We have the following commutative diagram with exact rows where the middle vertical map is an isomorphism by (1):
\[
\begin{tikzcd}
    0 \arrow{r} & \mathcal{E}^* \arrow{r} \arrow[d] & H^0(\mathcal{E})^*\otimes\mathcal{O}_{X} \arrow{r} \arrow[d,"\cong"] & M_{\mathcal{E}}^*\arrow{r}\arrow[d, equal] & 0\\
   0 \arrow{r} & M_{M_{\mathcal{E}}^*} \arrow{r} & H^0(M_{\mathcal{E}}^*)\otimes\mathcal{O}_{X} \arrow{r} &M_{\mathcal{E}}^*\arrow{r} & 0 
\end{tikzcd} \]
The assertion follows immediately from the above using Snake Lemma.\QEDB\par 

\begin{proposition}\label{ext1}
Assume $X$ is regular (i.e., $H^1(\mathcal{O}_X)=0)$, $H^1(\mathcal{E})=0$, and either (A) or (B) holds.
\begin{itemize}
    \item[(A)] If $n=2$, then $H^0(K_X)\otimes H^0(\mathcal{E})\to H^0(K_X\otimes \mathcal{E})$ surjects.
    \item[(B)] If $n\geq 3$, then $H^{n-2}(K_X\otimes \mathcal{E})=0$.
\end{itemize}
Then, the following holds.
\begin{itemize}
    \item[(1)] $H^1(M_{\mathcal{E}}^*)=0$.
    \item[(2)] If $H^{n-1}(K_X\otimes\mathcal{E})=0$ and $H^{n}(K_X\otimes\mathcal{E})=0$, then $\textrm{Ext}^1(M_{\mathcal{E}},M_{\mathcal{E}})\cong\textrm{Ext}^1(\mathcal{E},\mathcal{E})$.
\end{itemize}
\end{proposition}
\noindent\textit{Proof.} (1) The argument is similar to that of \cite[Proposition 2]{Cam}. We first show that if (A) or (B) holds, then $H^1(M_{\mathcal{E}}^*)=0$. Indeed, notice that by Serre duality, $H^1(M_{\mathcal{E}}^*)\cong H^{n-1}(K_X\otimes M_{\mathcal{E}})$ and the claim follows by tensoring \eqref{defining} by $K_X$ and taking cohomology as $H^{n-1}(K_X)\cong H^1(\mathcal{O}_X)=0$ by hypothesis. 

(2) Thanks to (1) and \eqref{referee}, we have $H^1(M_{\mathcal{E}}\otimes M_{\mathcal{E}}^*)\cong \textrm{Coker}(H^0(\xi))$. Moreover, once again by \eqref{referee}, we obtain
$\textrm{Coker}(H^0(\eta))\cong H^1(\mathcal{E}\otimes\mathcal{E}^*)$ since $H^1(\mathcal{E})=0$ by hypothesis. Thus the assertion follows by \eqref{cohreferee}.\QEDB\par 

\begin{remark}\label{irregular1}
If $X$ is arbitrary (i.e. not necessarily regular), assume $H^1(\mathcal{E})=0$, $H^{n-1}(K_X\otimes\mathcal{E})=0$, and $H^{n}(K_X\otimes\mathcal{E})=0$. Then the above proof shows $\textrm{Ext}^1(\mathcal{E},\mathcal{E})\subseteq \textrm{Ext}^1(M_{\mathcal{E}}, M_{\mathcal{E}})$. Indeed, it follows from \eqref{referee} that $\textrm{Coker}(H^0(\xi))\subseteq H^1(M_{\mathcal{E}}\otimes M_{\mathcal{E}}^*)$ and $\textrm{Coker}(H^0(\eta))\cong H^1(\mathcal{E}\otimes\mathcal{E}^*)$ whence the conclusion follows from \eqref{cohreferee}.
\end{remark}

\begin{proposition}\label{ext2}
Assume $X$ is regular, $H^1(\mathcal{E})=H^2(\mathcal{E})=0$, and either (A) or (B) holds.
\begin{itemize}
    \item[(A)] If $n=2$, then $p_g(X):=h^0(K_X)=0$ and $H^0(K_X\otimes\mathcal{E})=0$.
    \item[(B)] If $n\geq 3$, then $h^2(\mathcal{O}_X)=0$, and $H^{n-1}(K_X\otimes \mathcal{E})=0$.
\end{itemize}
Then $\textrm{Ext}^2(M_{\mathcal{E}},M_{\mathcal{E}})\cong\textrm{Ext}^2(\mathcal{E},\mathcal{E})$.
\end{proposition}
\noindent\textit{Proof.} It follows from Proposition \ref{ext1} that $H^1(M_{\mathcal{E}}^*)=0$. We aim to show that $H^2(M_{\mathcal{E}}^*)=0$. By Serre duality, $H^2(M_{\mathcal{E}}^*)\cong H^{n-2}(K_X\otimes M_{\mathcal{E}})$. To see the vanishing of this group, tensor \eqref{defining} by $K_X$ and take the long exact sequence of cohomology. For surfaces, the claim follows from the assumption $p_g(X)=0$. In higher dimension, the claim follows from the assumptions $H^{n-1}(K_X\otimes \mathcal{E})=0$ and $H^{n-2}(K_X)\cong H^2(\mathcal{O}_X)=0$.

Now, passing to the cohomology of the bottom exact sequence of \eqref{referee}, we obtain $H^2(M_{\mathcal{E}}\otimes M_{\mathcal{E}}^*)\cong H^1(\mathcal{E}\otimes M_{\mathcal{E}}^*)$. Moreover, the cohomology sequence associated with the right vertical exact sequence of \eqref{referee} yields $H^1(\mathcal{E}\otimes M_{\mathcal{E}}^*)\cong H^2(\mathcal{E}\otimes\mathcal{E}^*)$ thanks to the assumptions $H^i(\mathcal{E})=0$ for $i=1,2$. The assertion follows from combining these two isomorphisms.\QEDB\par  

\begin{remark}
If $X$ is arbitrary (i.e. not necessarily regular), assume $H^1(\mathcal{E})=0$ and $H^2(\mathcal{O}_X)=0$. Further, if $n\geq 3$, assume $\textrm{Ext}^2(\mathcal{E},\mathcal{E})=0$ and $H^{n-3}(K_X\otimes\mathcal{E})=0$. Then $\textrm{Ext}^2(M_{\mathcal{E}}, M_{\mathcal{E}})=0$.
Indeed, it follows from our assumption and the cohomology sequence associated to \eqref{defining} that $H^2(M_{\mathcal{E}})=0$. Consequently, the assertion follows for surfaces by the cohomology sequence associated to the left column of \eqref{referee}. 
In order to prove the assertion for $n\geq 3$, one needs to argue that $H^3(M_{\mathcal{E}}\otimes\mathcal{E}^*)=0$. But this follows via the cohomology sequence of the top row of \eqref{referee} thanks to our assumptions.
\end{remark}

Note that in the case of regular surfaces i.e. when $n=2$, we proved that $\textrm{Ext}^i(\mathcal{E},\mathcal{E})\cong\textrm{Ext}^i(M_{\mathcal{E}},M_{\mathcal{E}})$ for $i=0,1,2$ when $K_X\otimes\mathcal{E}$ has no cohomology. Thus, we are asking when $\mathcal{E}$ is stable that the moduli of stable bundles containing $K_X\otimes\mathcal{E}$ satisfies weak Brill--Noether which is a strong restriction. One way to ensure this is to require for embedded varieties that $K_X$ and the hyperplane bundle $H$ are related to each other and $\mathcal{E}$ satisfies ``nice'' conditions with respect to $H$. Keeping that in mind, we introduce a few definitions below.

\begin{definition}
Let $X\subseteq\mathbb{P}^N$ be a smooth projective variety and let $\mathcal{O}_X(1):=\mathcal{O}_{\mathbb{P}^N}(1)|_X$. Further, let $\mathcal{F}$ be a coherent sheaf on $X$.
\begin{itemize}
    \item[--] $X\subseteq\mathbb{P}^N$ is called {\it $s$--subcanonical} if $K_{X}=\mathcal{O}_{\mathbb{P}^N}(s)|_X$. When the embedding is projectively normal, it is called {\it arithmetically Cohen--Macaulay} (ACM for short) if $H^i(\mathcal{O}_X(j))=0$ for all $1\leq i\leq \dim X-1$ and for all $j\in\mathbb{Z}$. If $X\subseteq\mathbb{P}^N$ is ACM and $s$--subcanonical, then it is called {\it arithmetically Gorenstein} (AG for short).  
    \item[--] $\mathcal{F}$ is called {\it $m$--regular with respect to $\mathcal{O}_{X}(1)$} if $H^i(\mathcal{F}(m-i))=0$ for all $i\geq 1$.
    \item[--] $\mathcal{F}$ is called {\it normalized with respect to $\mathcal{O}_{X}(1)$} if $H^0(\mathcal{F}(-1))=0$ but $H^0(\mathcal{F})\neq 0$. 
\end{itemize}
\end{definition}

Recall the fundamental result of Mumford: any coherent sheaf $\mathcal{F}$ on $X\subseteq\mathbb{P}^N$ that is $0$--regular with respect to $\mathcal{O}_X(1)$ is globally generated. Further, if $\mathcal{F}$ is $m$--regular, then it is $(m+k)$--regular for all $k\geq 0$ (see for example \cite[Theorem 1.8.5]{Laz}). We will use these facts repeatedly, without any further reference. A key thing about $0$--regular bundles is the Castelnuovo--Mumford lemma that we will recall in Lemma \ref{cm}.

\begin{corollary}\label{subcan}
Let $X\subseteq\mathbb{P}^N$ be an $s$--subcanonical smooth regular projective variety of dimension $n$ with $H^2(\mathcal{O}_X)=0$ and set $H:=\mathcal{O}_{\mathbb{P}^N}(1)|_X$. Let $\mathcal{E}$ be a vector bundle on $X$ that is $0$--regular with respect to $H$. Further, assume that one of the following holds.
\begin{itemize}
    \item[(A)] $n=2$, $\mathcal{E}$ is normalized with respect to $H$, and $s= -1$.
    \item[(B)] $n\geq 3$ and $s\geq -(n-2)$.
\end{itemize}Then,
\begin{itemize}
    \item[(1)] $\mathcal{E}$ is globally generated and $\mathcal{E}\cong \left(M_{M_{\mathcal{E}}^*}\right)^*$.
    \item[(2)] $\textrm{Ext}^i(M_{\mathcal{E}},M_{\mathcal{E}})\cong\textrm{Ext}^i(\mathcal{E},\mathcal{E})$ for $i=0,1,2$.
\end{itemize}
\end{corollary}
\noindent\textit{Proof.} (1) We just discussed why $\mathcal{E}$ is globally generated. The fact that $\mathcal{E}\cong \left(M_{M_{\mathcal{E}}^*}\right)^*$ follows from Proposition \ref{hom} (3). Indeed, $\mathcal{E}$ is $m$--regular for all $m\geq 0$ and in particular $H^{n-1}(\mathcal{E}(l_1))=0$ for all $l_1\geq -(n-1)$ and $H^{n}(\mathcal{E}(l_2))=0$ for all $l_2\geq -n$. \par 
(2) As before $H^i(\mathcal{E})=0$ for all $i\geq 1$ due to the regularity assumption. Assume first that $n=2$ and $s=-1$. Then $p_g(X)=0$, $H^0(K_X\otimes\mathcal{E})=0$ since $\mathcal{E}$ is normalized by assumption and the assertion follows from Propositions \ref{hom} (2), \ref{ext1} and \ref{ext2}.\par 
Now if $n\geq 3$, regularity again implies the vanishings of $H^{n-1}(K_X\otimes\mathcal{E})$ and $H^{n-2}(K_X\otimes\mathcal{E})$ since $s\geq -(n-2)$ and the conclusion follows from Propositions \ref{hom} (2), \ref{ext1} and \ref{ext2}. That completes the proof.\QEDB\par 

\section{Semistability for (iterated--)syzygy bundles of Ulrich bundles}\label{3}

Unless otherwise stated, throughout this section, $X\subseteq\mathbb{P}^N$ is a smooth projective variety of degree $d$ and dimension $n$, $H:=\mathcal{O}_{\mathbb{P}^N}(1)|_X$, $\mathcal{E}$ is a vector bundle on $X$ of rank $r$. Further, let $X_n=X$ and let $X_{n-i}$ be the (smooth, irreducible) subvariety obtained by intersecting general sections $s_1,\cdots, s_i\in |H|$ for $0\leq i\leq n-1$. It follows from adjunction (and induction) that for all $0\leq i\leq n-1$
\begin{equation}\label{adjunction}
    K_{X_{n-i}}=(K_X+iH)|_{X_{n-i}}.
\end{equation}

\subsection{(Semi--)Stability of vector bundles} The {\it slope} of a torsion--free sheaf $\mathcal{F}$ with respect to $H$ is $$\mu_H(\mathcal{F}):=c_1(\mathcal{F})\cdot H^{n-1}/\textrm{rank}(\mathcal{F}).$$ 
Notice that $\mu_H(\mathcal{E})=\mu_{H|_{X_{n-i}}}(\mathcal{E}|_{X_{n-i}})$ for all $0\leq i\leq n-1$. Further, for any torsion--free sheaf $\mathcal{F}$, $P_{\mathcal{F}}$ denotes its Hilbert polynomial with respect to $H$.

\begin{definition} We define the notions of slope and Gieseker (semi--)stability.
\begin{itemize}
    \item[--] $\mathcal{E}$ is called {\it slope--stable (resp. slope--semistable) with respect to $H$} if for all subsheaf $\mathcal{F}\subsetneq \mathcal{E}$ with $0<\textrm{rank}(\mathcal{F})<\textrm{rank}(\mathcal{E})$, the inequality $\mu_H(\mathcal{F})<\mu_H(\mathcal{E})$ (resp. $\mu_H(\mathcal{F})\leq \mu_H(\mathcal{E})$) holds.
    \item[--] $\mathcal{E}$ is called {\it stable (resp. semistable) with respect to $H$} (in the sense of Gieseker) if for all subsheaf $\mathcal{F}\subsetneq \mathcal{E}$ with $0<\textrm{rank}(\mathcal{F})<\textrm{rank}(\mathcal{E})$, the inequality $P_{\mathcal{F}}/\textrm{rank}(\mathcal{F})<P_{\mathcal{E}}/\textrm{rank}(\mathcal{E})$ (resp. $P_{\mathcal{F}}/\textrm{rank}(\mathcal{F})\leq P_{\mathcal{E}}/\textrm{rank}(\mathcal{E})$) holds.
\end{itemize}
\end{definition}

On curves, (semi--)stability coincides with slope--(semi--)stability but in general we have the following chain of implications: 

\vspace{5pt}

\begin{center}
    slope--stable $\implies$ stable $\implies$ semistable $\implies$ slope--semistable.
\end{center}
\vspace{5pt}
Moreover, if $c_1(\mathcal{E})\cdot H^{n-1}$ and $\textrm{rank}(\mathcal{E})$ are coprime, the slope--semistability and slope--stability coincides.

\subsection{(Semi--)Stability for syzygy bundles of Ulrich bundles} We recall the definition of Ulrich bundles.

\begin{definition}
The vector bundle $\mathcal{E}$ is called {\it Ulrich with respect to $H$} if $H^j(\mathcal{E}(-j))=0$ for all $j>0$ (i.e., $\mathcal{E}$ is $0$--regular with respect to $H$) and $H^j(\mathcal{E}(-j-1))=0$ for all $j<n$.
\end{definition}

The definition above is equivalent to the definitions we gave in the introduction. We further recall some properties of Ulrich bundles which will often be used without any further reference.

\begin{proposition}\label{properties}
Assume $\mathcal{E}$ is Ulrich with respect to $H$. Then $h^0(\mathcal{E})=rd$. Further, for all $0\leq i\leq n-1$,
\begin{itemize}
    \item[(1)] $\mathcal{E}|_{X_{n-i}}$ is Ulrich with respect to $H|_{X_{n-i}}$ and globally generated,
    \item[(2)] $\mathcal{E}|_{X_{n-i}}$ is semistable with respect to $H|_{X_{n-i}}$,
    \item[(3)] $\mu_{H|_{X_{n-i}}}(\mathcal{E}|_{X_{n-i}})=d+g-1$ where $g$ is the genus of $X_1$. 
\end{itemize}
\end{proposition}

\noindent\textit{Proof.} The dimension of the space of global sections follows from \cite[Corollary 2.5 (ii)]{CKM}  (also from \cite[Lemma 2.4 (ii)]{CHGS}). (1) follows from \cite[Lemma 2.4 (i)]{CHGS},  induction, and \cite[Corollary 2.5 (i)]{CKM}. (2) follows from \cite[Theorem 2.9]{CHGS} (also from \cite[Proposition 2.6]{CKM}) and (1). Finally, (3) follows from \cite[Lemma 2.4 (iii)]{CHGS}  since $\mu_H(\mathcal{E})=\mu_{H|_{X_{n-i}}}(\mathcal{E}|_{X_{n-i}})$.\QEDB\par 

\vspace{5pt}

The following lemma shows that for Ulrich bundles, restriction of syzygy bundles are isomorphic to the syzygy bundles of the restrictions. In fact, being Ulrich is not required, weaker assumptions (i.e. normalized and $H^1(\mathcal{E}(-H))=0$) suffice. 

\begin{lemma}\label{meres}
Let $\mathcal{E}$ be Ulrich with respect to $H$ and assume $n\geq 2$. Then $M_{\mathcal{E}}|_{X_{i}}\cong M_{\mathcal{E}|_{X_{i}}}$ for all $1\leq i\leq n-1$.
\end{lemma}
\noindent\textit{Proof.} We have the following commutative diagram with exact rows.
\[
\begin{tikzcd}
    0 \arrow{r} & M_{\mathcal{E}}|_{X_{i}} \arrow{r} \arrow[d] & H^0(\mathcal{E})\otimes\mathcal{O}_{X_i} \arrow{r} \arrow{d} & \mathcal{E}|_{X_i}\arrow{r}\arrow[d, equal] & 0\\
   0 \arrow{r} & M_{\mathcal{E}|_{X_{i}}} \arrow{r} & H^0(\mathcal{E}|_{X_i})\otimes\mathcal{O}_{X_i} \arrow{r} &\mathcal{E}|_{X_i}\arrow{r} & 0 
\end{tikzcd} \]
Since Ulrich bundles are normalized, $0$--regular, and restriction of an Ulrich bundle to a hyperplane section remains Ulrich (see Proposition \ref{properties} (1)), it follows from induction that the middle vertical map is an isomorphism. The conclusion follows from Snake Lemma.\QEDB\par

\vspace{5pt}

A result of Butler gives a sufficient condition for semistability of syzygy bundles on curves:

\begin{theorem}\label{but}
(\cite[Theorem 1.2]{But94}) Let $\mathcal{F}$ be a semistable vector bundle on an irreducible smooth projective curve $C$ of genus $g$ with $\mu(\mathcal{F})\geq 2g$. Then $M_{\mathcal{F}}$ is semistable.
\end{theorem}

It is easy to see that if $\mathcal{F}$ is semistable on a smooth projective curve $C$ with $\mu(\mathcal{F})\geq 2g$ then $H^1(\mathcal{F})=0$ (see for e.g. \cite[Lemma 1.12 (2)]{But94}) and in this case $\mu(M_{\mathcal{F}})\geq -2$. Now we prove a result on semistability of syzygy bundles of Ulrich bundles.

\begin{proposition}\label{semistability}
Assume $\mathcal{E}$ is Ulrich with respect to $H$ and $n\geq 2$. If $(3-n)H^n>H^{n-1}K_X+2$ then $M_{\mathcal{E}}$ is slope--semistable with respect to $H$.
\end{proposition}

\noindent\textit{Proof.} It is enough to show that $M_{\mathcal{E}}|_{X_1}$ is semistable. By Lemma \ref{meres} $M_{\mathcal{E}}|_{X_1}\cong M_{\mathcal{E}|_{X_1}}$ and by Proposition \ref{properties} (2), $\mathcal{E}|_{X_1}$ is semistable on $X_1$. Thus, by Theorem \ref{but}, we need to show that $\mu(\mathcal{E}|_{X_1})\geq 2g(X_1)$. Thanks to Proposition \ref{properties} (3), it is enough to check that $d\geq g(X_1)+1$ i.e., $H^n\geq g(X_1)+1$. 

We know that $(H|_{X_2})^2+H|_{X_2}\cdot K_{X_2}=2g(X_1)-2$. Applying \eqref{adjunction} we obtain \begin{equation}\label{g(X_1)}
    H^n+H^{n-1}K_X+(n-2)H^n=2g(X_1)-2\implies g(X_1)=\frac{(n-1)H^{n}+H^{n-1}K_X}{2}+1.
\end{equation} Thus, it is enough to show that $H^n>\frac{(n-1)H^{n}+H^{n-1}K_X}{2}+1$ which follows from our assumption.\QEDB

\begin{corollary}\label{coprimes}
Let $L$ be an Ulrich line bundle with respect to $H$. Assume $n\geq 2$ and
\begin{itemize}
    \item[(A)] $(3-n)H^n>H^{n-1}K_X+2$,
    \item[(B)] $\textrm{gcd}(H^n-1,((n-1)H^{n}+H^{n-1}K_X)/2+1)=1$.
\end{itemize} 
Then $M_L$ slope--stable with respect to $H$.
\end{corollary}
\noindent\textit{Proof.} This follows immediately from Proposition \ref{semistability}. Indeed, if $L$ is an Ulrich line bundle with respect to $H$, then $$c_1(M_L)\cdot H=-L\cdot H=-\mu_H(L)=-(d+g(X_1)-1)$$ where $g(X_1)$ is the genus of $X_1$ and the last equality is obtained by Proposition \ref{properties} (3). On the other hand, $\textrm{rank}(M_L)=h^0(L)-1=d-1$ thanks to Proposition \ref{properties}. Thus, the assertion follows if $\textrm{gcd}(d+g(X_1)-1,d-1)$ which is the same as $\textrm{gcd}(d+g(X_1)-1,g(X_1))$ is $1$. Clearly, $d=H^n$ and the value for $g(X_1)$ is obtained from \eqref{g(X_1)}. The proof is now complete.\QEDB

\subsection{Koszul property and iterated syzygy bundles}\label{33} Next we define the iterated syzygy bundles $\mathcal{S}_k(\mathcal{E})$.

\begin{definition}
As before, let $\mathcal{E}$ be a vector bundle on $X$.
\begin{itemize}
    \item[--] We set $\mathcal{S}_{-1}(\mathcal{E}):=\mathcal{E}$.
    \item[--] If $\mathcal{E}$ is globally generated, set $\mathcal{S}_0(\mathcal{E}):=M_{\mathcal{E}}\otimes H$, otherwise $\mathcal{S}_0(\mathcal{E})$ does not exist.
    \item[--] We define $\mathcal{S}_{k}$ inductively: if for $k\geq 1$, $\mathcal{S}_{k-1}$ exists and is globally generated, define $\mathcal{S}_k(\mathcal{E}):=M_{\mathcal{S}_{k-1}(\mathcal{E})}\otimes H$, otherwise $\mathcal{S}_k(\mathcal{E})$ does not exist.
\end{itemize}
\end{definition}

The following lemma follows immediately from the above definition and \eqref{defining}.

\begin{lemma}\label{propertys}
Assume $\mathcal{S}_k(\mathcal{E})$ exists for some $k\geq 0$. Then for all $0\leq p\leq k$
\begin{itemize}
    \item[(1)] $\mathcal{S}_p(\mathcal{E})$ exists and $\mathcal{S}_{p'}(\mathcal{E})$ is globally generated for all $-1\leq p'\leq k-1$,
    \item[(2)] $\mathcal{S}_p(\mathcal{E})(-H)\cong M_{\mathcal{S}_{p-1}(\mathcal{E})}$ and one has the following short exact sequence 
    \begin{equation}\label{definings}
        0\to \mathcal{S}_p(\mathcal{E})(-H)\to H^0(\mathcal{S}_{p-1}(\mathcal{E}))\otimes\mathcal{O}_X\to \mathcal{S}_{p-1}(\mathcal{E})\to 0.
    \end{equation}
    \item[(3)]$H^0(\mathcal{S}_p(\mathcal{E})(-H))=0$. Moreover, if $H^1(\mathcal{O}_X)=0$, then $H^1(\mathcal{S}_p(\mathcal{E})(-H))=0$. 
\end{itemize}
\end{lemma}

Define the graded $\mathbb{C}$--algebra $R(H)$ as $R(H):=\oplus_{m\geq 0}H^0(mH)$. Recall when $X\subseteq\mathbb{P}^N$ is projectively normal (i.e., $H^0(\mathcal{O}_{\mathbb{P}^N}(m))\to H^0(mH)$ surjects for all $m$), $R(H)$ is the homogeneous coordinate ring of the embedding. The following characterization of $\Gamma_*(\mathcal{E})$ being Koszul follows directly by the argument of Lazarsfeld (see \cite{Par}). We provide the proof for the sake of completeness.

\begin{proposition}\label{critkos}
Let $X\subseteq\mathbb{P}^N$ be projectively normal with homogeneous coordinate ring $R(H)$. Assume $\mathcal{E}$ is a vector bundle for which $\mathcal{S}_k(\mathcal{E})$ exists for all $k\geq 0$. Then $\Gamma_*(\mathcal{E})$ is a Koszul $R(H)$ module if and only if $H^0(\mathcal{S}_{k}(\mathcal{E}))\otimes H^0(mH)\to H^0(\mathcal{S}_{k}(\mathcal{E})(mH))$ surjects for all $k\geq -1, m\geq 0$. 
\end{proposition}

\noindent\textit{Proof.} Recall from the Introduction that $M:=\Gamma_*(\mathcal{E})=\oplus_{m\geq 0}M_m$ is a Koszul $R:=R(H)=\oplus_{m\geq 0} R_m$ module if and only if $M$ admits a linear minimal graded free resolution $$\cdots\to E_i\to\cdots \to E_1\to E_0\to M\to 0.$$
Observe that $E_0=\oplus R$ if and only if $R_m\otimes M_0\to M_m$ surjects for all $m\geq 0$ which is equivalent to the surjectivity of $H^0(\mathcal{E})\otimes H^0(mH)\to H^0(\mathcal{E}(mH))$ for all $m\geq 0$. Define $R_m^1:=\textrm{Ker}(R_m\otimes M_0\to M_m)$ and $R^1:=\oplus_{m\geq 0}R_m^1$. Then $E_1=\oplus R(-1)$ if and only if $R^1$ is generated in degree one i.e., $R_m\otimes R^1_1\to R^1_{m+1}$ surjects for all $m\geq 0$. Notice that $R^1_1=H^0(M_{\mathcal{E}}(H))$ and consequently the surjection $R_m\otimes R^1_1\to R^1_{m+1}$ is equivalent to the surjection $H^0(M_{\mathcal{E}}(H))\otimes H^0(mH)\to H^0(M_{\mathcal{E}}(m+1)H)$. An easy induction finishes the proof.\QEDB\par  

\vspace{5pt}

We now define the property $(\textrm{Kos}_k)$ for vector bundles $\mathcal{E}$ that is crucial for the proof of Theorem \ref{main}.

\begin{definition}
Let $k\geq 0$ be an integer. We say that $\mathcal{E}$ satisfies $(\textrm{Kos}_k)$ if the following two condition hold
\begin{itemize}
    \item[(1)] $\mathcal{S}_p(\mathcal{E})$ exists for all $-1\leq p\leq k-1$, and 
    \item[(2)] $H^0(\mathcal{S}_p(\mathcal{E}))\otimes H^0(H)\to H^0(\mathcal{S}_p(\mathcal{E})\otimes H)$ surjects for all $-1\leq p\leq k-1$.
\end{itemize}
If $\mathcal{E}$ satisfies $(\textrm{Kos}_k)$ for all $k\geq0$, we say that $\mathcal{E}$ satisfies $(\textrm{Kos}_{\infty})$
\end{definition}

In what follows, we will require to prove surjectivity of multiplication maps. The following observation of Gallego and the first author will be useful for us.

\begin{lemma}\label{gpobs}
(\cite[Observation 1.4.1]{GP99}) Let $\mathcal{F}$ and $L_0:=\mathcal{O}_X$, $L_1$, $L_2,\cdots,L_r$ be coherent sheaves on a variety $X$. Assume the following evaluation maps $$H^0\left(\mathcal{F}\otimes \bigotimes\limits_{i=0}^{k-1}L_i\right)\otimes H^0(L_k)\to H^0\left(\mathcal{F}\otimes \bigotimes\limits_{i=1}^{k}L_i\right)$$ surjects for all $1\leq k\leq r$. Then the evaluation map  $H^0(\mathcal{F})\otimes H^0\left(\bigotimes\limits_{i=1}^{r}L_i\right)\to H^0\left(\mathcal{F}\otimes \bigotimes\limits_{i=1}^{r}L_i\right)$ also surjects.
\end{lemma}

We also recall the following lemma known as Castelnuovo--Mumford lemma that provides a very helpful criterion for surjectivity of multiplication maps for $0$--regular bundles that we alluded to earlier.

\begin{lemma}\label{cm}
(\cite{Mu}) Let $L$ be a base point free line bundle on a variety $X$ and let $\mathscr{F}$ be a coherent sheaf on $X$. If $H^i(\mathscr{F}\otimes L^{-i})=0$ for all $i\geq 1$ (i.e., $\mathcal{F}$ is $0$--regular with respect to $L$) then the multiplication map $H^0(\mathscr{F}\otimes L^{\otimes i})\otimes H^0(L)\rightarrow H^0(\mathscr{F}\otimes L^{\otimes i+1})$ surjects for all $i\geq 0$.
\end{lemma} 

An immediate consequence of the following lemma is the $0$--regularity of $\mathcal{S}_p(\mathcal{E})(H)$ with respect to $H$.

\begin{lemma}\label{propertykos}
Assume $\mathcal{O}_X$ is $2$--regular and $\mathcal{E}$ is $1$--regular with respect to $H$. Further assume $\mathcal{E}$ satisfies $(\textrm{Kos}_k)$ for some $k\geq 0$. Then $\mathcal{S}_p(\mathcal{E})$ is $1$--regular with respect to $H$ for all $-1\leq p\leq k-1$.
\end{lemma}
\noindent\textit{Proof.} There is nothing to prove if $p=-1$, so we assume $p\geq 0$. We use induction on $i$ to show that $H^i(\mathcal{S}_p(\mathcal{E})(-(i-1)H))=0$ when $0\leq p\leq k-1$. This vanishing for $i=1$ can be seen by tensoring \eqref{definings} by $H$ and using the definition of property $(\textrm{Kos}_k)$ and the regularity assumption on $\mathcal{O}_X$. Similarly, for higher $i$ the vanishing follows inductively by twisting \eqref{definings} by $-(i-2)H$ and using the regularity assumption of $\mathcal{O}_X$.\QEDB

\begin{corollary}\label{kosiff2}
Assume $\mathcal{O}_X$ is $2$--regular and $\mathcal{E}$ is $1$--regular with respect to $H$. Further assume $\mathcal{S}_k(\mathcal{E})$ exists for all $k\geq 0$. Then $\mathcal{E}$ satisfies property $(\textrm{Kos}_{\infty})$ if and only if $\Gamma_*(\mathcal{E})$ is a Koszul $R(H)$ module.
\end{corollary}

\noindent\textit{Proof.} Notice that $\mathcal{S}_k(\mathcal{E})(H)$ is $0$--regular with respect to $H$ by Lemma \ref{propertykos}. Thus, if $\mathcal{E}$ satisfies $(\textrm{Kos}_{\infty})$ then $H^0(\mathcal{S}_k(\mathcal{E}))\otimes H^0(mH)\to H^0(\mathcal{S}_k(\mathcal{E})(mH))$ surjects by Lemma \ref{gpobs} and Lemma \ref{cm}. Consequently $\Gamma_*(\mathcal{E})$ is Koszul $R(H)$ module by Proposition \ref{critkos}. The converse follows directly from Proposition \ref{critkos}.\QEDB

\begin{proposition}\label{restricts}
Let $X$ be regular and let $\mathcal{E}$ be an Ulrich vector bundle with respect to $H$. Assume $\mathcal{S}_k(\mathcal{E})$ exists for some $k\geq 0$. Then for all $0\leq p\leq k$ and for all $1\leq i\leq n-1$,
 $\mathcal{S}_p(\mathcal{E}|_{X_i})$ exists and 
 $\mathcal{S}_p(\mathcal{E}|_{X_i})\cong \mathcal{S}_p(\mathcal{E})|_{X_i}$.
\end{proposition}
\noindent\textit{Proof.} First assume $p=0$. Then, by hypothesis $\mathcal{E}$ is globally generated and consequently $\mathcal{E}|_{X_i}$ is globally generated. Thus, $\mathcal{S}_0(\mathcal{E}|_{X_i})$ exists. Further $\mathcal{S}_0(\mathcal{E}|_{X_i})\cong \mathcal{S}_0(\mathcal{E})|_{X_i}$ follows from Lemma \ref{meres}.

Now, assume the statements hold for all $0\leq p\leq p_0\leq k-1$. We know that $\mathcal{S}_{p}(\mathcal{E})$ exists and globally generated for all $p\leq p_0$ by Lemma \ref{propertys} (1). Thus $\mathcal{S}_{p_0}(\mathcal{E}|_{X_i})\cong \mathcal{S}_{p_0}(\mathcal{E})|_{X_i}$ by induction hypothesis and hence globally generated $\implies \mathcal{S}_{p_0+1}(\mathcal{E}|_{X_i})$ exists. By Lemma \ref{propertys} (2) we have the following commutative diagram with exact rows.
\[
\begin{tikzcd}
    0 \arrow{r} & \mathcal{S}_{p_0+1}(\mathcal{E})|_{X_{i}} \arrow{r} \arrow[d] & H^0(\mathcal{S}_{p_0}(\mathcal{E}))\otimes H|_{X_i} \arrow{r} \arrow{d} & \mathcal{S}_{p_0}(\mathcal{E})(H)|_{X_i}\arrow{r}\arrow[d, equal] & 0\\
   0 \arrow{r} & \mathcal{S}_{p_0+1}(\mathcal{E}|_{X_{i}}) \arrow{r} & H^0(\mathcal{S}_{p_0}(\mathcal{E}|_{X_i}))\otimes H|_{X_i} \arrow{r} &\mathcal{S}_{p_0}(\mathcal{E})(H)|_{X_i}\arrow{r} & 0 
\end{tikzcd} \]
We claim that the middle vertical map is an isomorphism. Indeed, thanks to our induction hypothesis, it is enough  to show that $H^0(\mathcal{S}_{p_0}(\mathcal{E}|_{X_j})(-H|_{X_j}))=H^1(\mathcal{S}_{p_0}(\mathcal{E}|_{X_j})(-H|_{X_j}))=0$ for all $i+1\leq j\leq n$. By Kodaira Vanishing Theorem, each $X_j$ is regular and consequently the vanishings follow from Lemma \ref{propertys} (3).\QEDB

\vspace{5pt}

By repeated application of the following result of Gallego and the first author, we will boil the problem of surjectivity of multiplication map on a variety down to a curve section.

\begin{lemma}\label{gpreg}
 (\cite[Observation 2.3]{GP99}) Let $X$ be a regular variety (i.e. $H^1(\mathscr{O}_X)=0$). Let $E$ be a vector bundle and let $D$ be a divisor such that $L=\mathscr{O}_X(D)$ is globally generated and $H^1(E\otimes L^*)=0$. If the multiplication map $H^0(E\vert_D)\otimes H^0(L\vert_D)\rightarrow H^0((E\otimes L)\vert_D)$ surjects then $H^0(E)\otimes H^0(L)\rightarrow H^0(E\otimes L)$ also surjects.
\end{lemma} 

Finally, when we have a multiplication map on a curve, we will use the following result of Butler.

\begin{theorem}\label{butler}
\cite[Theorem 2.1]{But94}) Let $E$ and $F$ be semistable vector bundles over a smooth projective curve $C$ of genus $g$ such that $\mu(E)\geq 2g$, and $\mu(F)> 2g$.
Then the multiplication map of global sections $H^0(E)\otimes H^0(F)\rightarrow H^0(E\otimes F)$ surjects.
\end{theorem} 

We use the above results to prove the following theorem the hypotheses of which for dimension $n\geq 3$ is satisfied for example when $X$ is a quadric hypersurface (Fano of index $n$ i.e. $(-n)$--subcanonical variety).

\begin{theorem}\label{mains}
Assume $n\geq 2$ and $X\subseteq\mathbb{P}^N$ is ACM. Let $\mathcal{O}_X$ be $2$--regular and $\mathcal{E}$ be Ulrich with respect to $H$. Further assume $(2-n)H^n\geq H^{n-1}K_X+4$. Then for all $p\geq -1$
\begin{itemize}
    \item[(1)] $\mathcal{S}_p(\mathcal{E})$ exists,
    \item[(2)] $\mu_{H|_{X_i}}(\mathcal{S}_p(\mathcal{E})|_{X_i})\geq 2g$ and $\mathcal{S}_p(\mathcal{E}|_{X_i})$ is slope--semistable with respect to $H|_{X_i}$ for all $1\leq i\leq n$,
    \item[(3)] $H^0(\mathcal{S}_p(\mathcal{E}|_{X_i}))\otimes H^0(H|_{X_i})\to H^0(\mathcal{S}_p(\mathcal{E}|_{X_i})\otimes H|_{X_i})$ surjects for all $1\leq i\leq n$.
\end{itemize}
In particular, $\mathcal{E}|_{X_i}$ satisfies $(\textrm{Kos}_{\infty})$ for all $1\leq i\leq n$.
\end{theorem}
\noindent\textit{Proof.} It follows from adjunction that $H|_{X_2}^2+H|_{X_2}K_{X_2}=2g-2$ and consequently by \eqref{adjunction}, it is easy to see that the given inequality is equivalent to $H^n\geq 2g+2$. We proceed by induction on $p$. 

First we settle the case for $p=-1$.  The assertion (1) is clear, (2) follows from Proposition \ref{properties} since $\mu_{H|_{X_i}}(\mathcal{E}|_{X_i})=\mu_H(\mathcal{E})\geq 3g+1$ by our assumption. Finally (3) follows from Proposition \ref{properties} and Lemma \ref{cm}.

Now assume $p=0$. Clearly $\mathcal{S}_0(\mathcal{E})=M_{\mathcal{E}}(H)$ exists. To prove (2), it is enough to prove it for $i=1$. Notice that $M_{\mathcal{E}}(H)|_{X_1}\cong M_{\mathcal{E}|_{X_1}}(H|_{X_1})$. Since $\mu(\mathcal{E}|_{X_1})\geq 2g$ and $\mathcal{E}|_{X_1}$ is Ulrich by Proposition \ref{properties} and thus semistable, it follows by Theorem \ref{but} that $M_{\mathcal{E}|_{X_1}}$ is semistable with 
$\mu(M_{\mathcal{E}|_{X_1}})\geq -2$ and thus $M_{\mathcal{E}|_{X_1}}(H|_{X_1})$ is semistable with $\mu(M_{\mathcal{E}|_{X_1}}(H|_{X_1}))\geq 2g$. To prove (3), notice that $X_i$ is regular for all $i\geq 2$ by Kodaira Vanishing Theorem. Thus by Lemma \ref{elt} (1) and Lemma \ref{gpreg}, it is enough to consider the case when $i=1$. But in this case, the required surjection follows from Theorem \ref{butler}.

Now assume (1), (2) and (3) hold for some $0\leq p\leq p_0$, we want to show that they hold for $p=p_0+1$. Consider the bundle $\mathcal{S}_{p_0}(\mathcal{E})(H)$, it follows from Lemma 
\ref{propertykos} 
that it is $0$--regular with respect to $H$ and consequently globally generated. Now since (3) holds with $p=p_0$, it follows that $\mathcal{S}_{p_0}(\mathcal{E})$ is globally generated and consequently $\mathcal{S}_{p_0+1}(\mathcal{E})$ exists. For (2), thanks to Proposition \ref{restricts}, it is enough to show that $\mathcal{S}_{p_0+1}(\mathcal{E}|_{X_1})$ is slope--semistable with $\mu(\mathcal{S}_{p_0+1}(\mathcal{E}|_{X_1}))\geq 2g$. By Lemma \ref{propertys} (2) $\mathcal{S}_{p_0+1}(\mathcal{E}|_{X_1})\cong M_{\mathcal{S}_{p_0}(\mathcal{E}|_{X_1})}(H|_{X_1})$. By hypothesis $\mathcal{S}_{p_0}(\mathcal{E}|_{X_1})$ is semistable with $\mu(\mathcal{S}_{p_0}(\mathcal{E}|_{X_1}))\geq 2g$. Consequently, by Theorem \ref{but}, $M_{\mathcal{S}_{p_0}(\mathcal{E}|_{X_1})}$ is semistable with slope $\geq -2$ and the conclusion follows since $\mu(H|_{X_1})\geq 2g+2$. Finally (3) follows by Lemma \ref{gpreg}, Theorem \ref{butler} and Lemma \ref{propertys} (3) since $X_i$ is regular for $i\geq 2$. \QEDB

\vspace{5pt}

It is interesting to note that for ACM varieties, $\Gamma_*(\mathcal{E})$ being a Koszul $R(H)$ module is equivalent to $R(H)$ itself being Koszul. We give an algebraic proof below. In the following proof, we will use the fact that a finitely generated graded module $M$ over a graded $\mathbb{C}$--algebra $R=\oplus _{m\geq 0}R_m$ with $R_0=\mathbb{C}$ and $R$ generated in degree one being Koszul is equivalent to $\textrm{Tor}^R_i(M,\mathbb{C})$ being purely of degree $i$.

\begin{theorem}\label{kosiff}
Let $\mathcal{E}$ be Ulrich with respect to $H$ on a projectively normal smooth ACM projective variety $X\subseteq\mathbb{P}^N$. Then $\Gamma_*(\mathcal{E})$ is Koszul with respect to the homogeneous coordinate ring $R(H)$ of the embedding $X\subseteq\mathbb{P}^N$ if and only if $R(H)$ is Koszul.
\end{theorem}

\noindent\textit{Proof.} Set $M=\Gamma_*(\mathcal{E})$ and recall that $\mathcal{E}$ being Ulrich is equivalent to $M$ being a linear maximal Cohen--Macaulay $S:=\mathbb{C}[X_0,\cdots,X_N]$ module. Now $R:=R(H)$ is Cohen--Macaulay by hypothesis. Also observe that since $M$ is a projective $R$ module, by Auslander--Buchsbaum Theorem, $M$ has maximal depth as $R$ module. Let $\mathfrak{m}$ be the irrelevant maximal ideal of $R$. We can find a regular sequence $x_1,\cdots , x_n$ of linear forms by prime avoidance i.e. choose $x_1\in R-(\mathfrak{p}_1\cup\cdots\cup \mathfrak{p}_r\cup\mathfrak{m}^2)$ ($\mathfrak{p}_i$'s are associated primes of $M$) and inductively choose the subsequent $x_i$'s. Since $M$ has maximal depth i.e. $\textrm{depth}(M)=\textrm{depth}(R)$, we have that $M/(x_1,\cdots,x_n)M$ is an Artinian module and since $x_i$'s are linear forms, $M/(x_1,\cdots,x_n)M=\oplus\mathbb{C}$ where we identify $\mathbb{C}$ as $R/\mathfrak{m}$. Thus, $M$ is Koszul $\iff$ $M/(x_1,\cdots,x_n)M$ is Koszul $\iff$ $\textrm{Tor}_i^R(M/(x_1,\cdots,x_n),\mathbb{C})$ is of pure degree $i$ $\iff$ $\oplus\textrm{Tor}_i^R(\mathbb{C},\mathbb{C})$ is of pure degree $i$ $\iff$ $\textrm{Tor}_i^R(\mathbb{C},\mathbb{C})$ is of pure degree $i$ $\iff$ $R$ is Koszul. That completes the proof.\QEDB

\begin{remark}\label{ship}
A related result in commutative algebra was proven by Iyengar and R\"omer in \cite{IR}, Theorem 3.4. Also, by the above theorem, under the given hypotheses, Koszulity of $\Gamma_*(\mathcal{E})$ for an Ulrich bundle $\mathcal{E}$ is completely determined by the Koszulity of the embedded variety. In \cite{KMS}, Kulkarni, Mustopa and Shipman introduced the notion of $\delta$--Ulrich sheaves: a reflexive sheaf $\mathcal{F}$ on a polarized variety $(X,H)$ is {\it $\delta$--Ulrich} if there exists a smooth one--dimensional linear section $Y$ of $X$ such that the restriction $\mathcal{F}|_Y$ is an Ulrich sheaf on $Y$. Even though it is not known whether every ACM variety admits an Ulrich sheaf, the authors in \cite{KMS} showed that any embedded normal ACM variety admits a $\delta$--Ulrich sheaf.
\end{remark}

\begin{remark}\label{rmkc}
By Theorems \ref{kosiff}, \ref{mains} and Corollary \ref{kosiff2}, we conclude that anticanonically embedded Del Pezzo surfaces $X_d\subseteq\mathbb{P}^d$ (see Subsection \ref{43}) are Koszul for $d\geq 4$ (which is also well--known). We remark that Theorem \ref{mains} shows the existence and slope--semistability of the iterated syzygy bundles. On the other hand, the above corollary also shows that for Ulrich bundles $\mathcal{E}$ on cubic Del Pezzo surfaces, $\Gamma_*(\mathcal{E})$ can not be a Koszul $R(H)$ module since if it is, then by Theorem \ref{kosiff} $R(H)$ will be Koszul. This is a contradiction since a cubic Del Pezzo surface is not normally presented. 
\end{remark}

The following corollary is the final result of this section.

\begin{corollary}\label{extis}
Let $X\subseteq\mathbb{P}^N$ be a smooth regular $(-1)$--subcanonical surface and let $H:=\mathcal{O}_{\mathbb{P}^N}(1)|_X$. Assume $-K_XH\geq 4$. Let $\mathcal{E}$ be an Ulrich vector bundle with respect to $H$. Then 
\begin{itemize}
    \item[(1)] $\mathcal{E}$ satisfies $(\textrm{Kos}_{\infty})$ and $\mathcal{S}_p(\mathcal{E})$ is slope--semistable for all $p\geq 1$,
    \item[(2)] $\mathcal{S}_p(\mathcal{E})\cong \left(M_{(\mathcal{S}_{p+1}(\mathcal{E}))^*(H)}\right)^*$ for all $p\geq -1$,
    \item[(3)] $\textrm{Ext}^i(\mathcal{S}_p(\mathcal{E}),\mathcal{S}_p(\mathcal{E}))\cong\textrm{Ext}^i(\mathcal{E},\mathcal{E})$ for $i=0,1,2$ and for all $p\geq -1$.
\end{itemize}
\end{corollary}
\noindent\textit{Proof.} (1) is an immediate consequence of Theorem \ref{mains}. The assertion (2) follows from Corollary \ref{subcan} (1) for $p=-1$. For $p\geq 0$, $H^1(\mathcal{S}_p(\mathcal{E})(-H))=H^2(\mathcal{S}_p(\mathcal{E})(-H))=0$ by Lemma \ref{propertys} (3) and Lemma \ref{propertykos}. The assertion follows from Lemma \ref{propertys} (2) and Proposition \ref{hom} (3). Now we prove (3). The assertion is immediate for $p=-1$ by Corollary \ref{subcan} (2). Assume it holds for $0\leq p\leq p_0$ and we want to show it for $p=p_0+1$. Notice that $\textrm{Ext}^i(\mathcal{S}_{p_0+1}(\mathcal{E}),\mathcal{S}_{p_0+1}(\mathcal{E}))=\textrm{Ext}^i(M_{\mathcal{S}_{p_0}(\mathcal{E})}, M_{\mathcal{S}_{p_0}(\mathcal{E})})$ by Lemma \ref{propertys} (2). Notice that $H^i({\mathcal{S}_{p_0}(\mathcal{E})}(-H))=0$ for $i=0,1,2$ by Lemma \ref{propertys} (3) and Lemma \ref{propertykos}, $H^i({\mathcal{S}_{p_0}(\mathcal{E})})=0$ for $i=1,2$ by Lemma \ref{propertykos} (the regularity condition of $\mathcal{O}_X$ follows from Kodaira Vanishing Theorem). It follows from Proposition \ref{hom}, Proposition \ref{ext1} and Proposition \ref{ext2} that $\textrm{Ext}^i(M_{\mathcal{S}_{p_0}(\mathcal{E})}, M_{\mathcal{S}_{p_0}(\mathcal{E})})=\textrm{Ext}^i({\mathcal{S}_{p_0}(\mathcal{E})}, {\mathcal{S}_{p_0}(\mathcal{E})})$ for $i=0,1,2$. The conclusion follows from the induction hypothesis.\QEDB

\section{Consequences on the moduli of vector bundles on Del Pezzo surfaces}\label{4}

The objective for this section is to prove the results stated in the introduction. 

\subsection{Chern classes of vector bundles} Let $Y$ be a smooth projective surface. We recall a few well--known formulae for Chern classes that will enable us to carry out our computations.

\subsubsection{} Let $\mathcal{F}, \mathcal{G}$ be vector bundles on $Y$ of rank $s$ and $t$ respectively. The first Chern class of their tensor product can be calculated using the following formula
\begin{equation}\label{c1tensor}
    c_1(\mathcal{F}\otimes\mathcal{G})=tc_1(\mathcal{F})+sc_1(\mathcal{G}).
\end{equation}
Also, the second Chern class of the tensor product can be obtained by the following formula
\begin{equation}\label{c2tensorv}
    c_2(\mathcal{F}\otimes\mathcal{G})=\frac{s(s-1)}{2}c_1(\mathcal{G})^2+sc_2(\mathcal{G})+(st-1)c_1(\mathcal{F})c_1(\mathcal{G})+tc_2(\mathcal{F})+\frac{t(t-1)}{2}c_1(\mathcal{F})^2.
\end{equation}

\subsubsection{} Let $\mathcal{F}_1,\mathcal{F}_2,\cdots,\mathcal{F}_r$ be vector bundles on $Y$. The first Chern class of their direct sum is given by 
\begin{equation}\label{c1sum}
    c_1\left(\bigoplus\limits_{i=1}^r\mathcal{F}_i\right)=\sum_{i=1}^rc_1(\mathcal{F}_i).
\end{equation}
We fix the notation $I_r:=\left\{1,\cdots,r\right\}$ for the rest of this section. The second Chern class of the direct sum can be calculated by the following formula 
\begin{equation}\label{c2sum}
    c_2\left(\bigoplus\limits_{i=1}^r\mathcal{F}_i\right)=\sum_{i=1}^rc_2(\mathcal{F}_i)+\sum_{i,j\in I_r,i<j}c_1(\mathcal{F}_i) c_1(\mathcal{F}_j).
\end{equation}

\subsubsection{} The Riemann--Roch formula for a vector bundle $\mathcal{F}$ of rank $s$ on $Y$ reads as follows
\begin{equation}\label{rr}
    \chi(\mathcal{F})=s\chi(\mathcal{O}_X)+\frac{1}{2}\left(c_1(\mathcal{F})^2-c_1(\mathcal{F})\cdot K_X\right)-c_2(\mathcal{F}).
\end{equation}

\subsubsection{} The discriminant of a vector bundle $\mathcal{F}$ of rank $s$ on $Y$ is by definition $$\Delta(\mathcal{F}):=2sc_2(\mathcal{F})-(s-1)c_1(\mathcal{F})^2.$$

\subsection{Deformations and moduli of vector bundles}\label{42} Let $Y$ be a smooth projective surface and let $\mathcal{F}$ be a vector bundle on $Y$. The infinitesimal deformations of $\mathcal{F}$ are parametrized by elements of $\textrm{Ext}^1(\mathcal{F},\mathcal{F})$. We denote by $\textrm{Def}(\mathcal{F})$ the versal deformation space of $\mathcal{F}$. $\mathcal{F}$ is called {\it unobstructed} if $\textrm{Def}(\mathcal{F})$ is smooth which is equivalent to $\dim\textrm{Def}(\mathcal{F})=\dim \textrm{Ext}^1(\mathcal{F},\mathcal{F})$. To check unobstructedness of $\mathcal{F}$, it is enough to show that $\textrm{Ext}^2(\mathcal{F},\mathcal{F})=0$. Recall that a torsion--free coherent sheaf $\mathcal{G}$on $Y$ is called {\it simple} if $\textrm{Hom}(\mathcal{G},\mathcal{G})\cong\mathbb{C}$. If $\mathcal{F}$ is a simple vector bundle of rank $s$ with $c_i(\mathcal{F})=c_i$ on a rational surface $Y$, and satisfies $\textrm{Ext}^2(\mathcal{F},\mathcal{F})=0$, then from Riemann--Roch it follows that
\begin{equation}
    \dim\textrm{Def}(\mathcal{F})=2sc_2-(s-1)c_1^2-(s^2-1)
\end{equation}
\begin{definition}
Let $D$ be an effective divisor on a smooth projective surface $Y$ and let $\mathcal{G}$ be a torsion--free coherent sheaf on $Y$. $\mathcal{G}$ is called {\it $D$--prioritary} if $\textrm{Ext}^2(\mathcal{G},\mathcal{G}(-D))=0$. For ${\bf v}=(s,c_1,c_2)$, we denote by 
\begin{itemize}
    \item[--] $\mathfrak{P}_{D}({\bf v})$ (or $\mathfrak{P}_{D}(s;c_1,c_2)$) the stack of $D$-prioritary torsion--free coherent sheaves on $Y$ of rank $s$ and Chern classes $c_i$,
    \item[--] $\mathfrak{Spl}_D({\bf v})$ (or $\mathfrak{Spl}_D(s;c_1,c_2)$) the moduli space of simple $D$--prioritary torsion--free coherent sheaves of rank $s$ and Chern classes $c_i$.
\end{itemize}  
\end{definition}

We remark that $\mathfrak{Spl}_D({\bf v})$ is an open substack of $\mathfrak{P}_{D}({\bf v})$.
Now fix an ample line bundle $H$ on $Y$ and let $\mathfrak{M}_H({\bf v})$ (or $\mathfrak{M}_H(s;c_1,c_2)$) be the coarse moduli space of slope--stable vector bundles on $Y$ with respect to $H$. It is easy to see that if $H$ satisfies $H\cdot (K_X+D)<0$ where $D$ is an effective divisor, then any slope--semistable vector bundle $\mathcal{F}$ is $D$--prioritary and satisfies $\textrm{Ext}^2(\mathcal{F},\mathcal{F})=0$.  Moreover, in this case, the moduli space $\mathfrak{M}_H({\bf v})$ is an open substack of $\mathfrak{P}_{D}({\bf v})$. The following fundamental theorem is due to Walter.

\begin{theorem}\label{wal}
(\cite{Wal}) Let $Y$ be a birationally ruled surface with fiber class
$F$ for a fixed ruling. Then the stack $\mathfrak{P}_{F} ({\bf v})$
is irreducible whenever it is non--empty. Moreover, if $s\geq 2$ where ${\bf v}=(s,c_1,c_2)$, then the general
element of $\mathfrak{P}_{F} ({\bf v})$ is a vector bundle. In particular, if $H$ is a polarization such that $H\cdot (K_X + F)< 0$ and $\mathfrak{M}_{H} ({\bf v})$ is non--empty, then $\mathfrak{M}_{H} ({\bf v})$ is smooth and irreducible. 
\end{theorem}

\smallskip

We will also require the notion of {\it modular families} that we define below following \cite{CHGS}.

\begin{definition}
A {\it modular family} for a class of vector bundles on a smooth projective variety $X$ is a flat family of vector bundles $\mathbb{E}$ on $X\times S/S$ where $S$ is a scheme of finite type such that 
\begin{itemize}
    \item[(1)] each isomorphism class of vector bundles occurs at least once and at most finitely many times;
    \item[(2)] for each $s\in S$, $\hat{\mathcal{O}}_{S,s}$ with the induced family pro--represents the local deformation functor;
    \item[(3)] for any other flat family $\mathbb{E}'$ on $X\times S'/S'$ there exists an \'etale surjective map $S''\to S'$ an a morphism $S''\to S$ such that $\mathbb{E}'\times _{S'}S''/X\times S''\cong\mathbb{E}\times _{S'}S''/X\times S''$.
\end{itemize}
\end{definition}

It follows from \cite[Proposition 2.10]{CHGS} that on a smooth projective variety $X$, any bounded family of simple bundles $\mathcal{E}$ with given rank and Chern classes satisfying $H^2
(\mathcal{E}\otimes \mathcal{E}^*) = 0$ has a smooth modular family.

\subsection{Moduli of vector bundles on Del Pezzo surfaces}\label{43} A Del Pezzo surface $X_d$ of degree $3\leq d\leq 9$ is a blow--up of $\mathbb{P}^2$ in $9-d$ points in general position. The anti--canonical bundle $-K_{X_d}$ is very ample if $d\geq 3$ and $h^0(-K_{X_d})=d+1$. In this section, we shall only consider anti--canonically embedded Del Pezzo surfaces $X_d\subseteq\mathbb{P}^d$. The Picard group of $X_d$ is $\mathbb{Z}^{\oplus 10-d}$ generated by $L, E_1,\cdots, E_{9-d}$ where $L$ is the pull--back of hyperplane class of $\mathbb{P}^2$ and $E_i$'s are exceptional divisors. The intersections are given by $$L^2=1,\quad L\cdot E_i=0,\quad E_i\cdot E_j=-\delta_{ij}.$$
The class of the canonical divisor is given by $$K_{X_d}=-3L+\sum_{i=1}^{9-d}E_i.$$ In what follows, for a Del Pezzo surface $X_d\subseteq\mathbb{P}^d$ with $d\leq 8$, we set $F:=L-E_1$; the complete linear series $|F|$ gives $X_d$ a structure of birationally ruled surface. We also set $H:=\mathcal{O}_{X_d}(1):=\mathcal{O}_{\mathbb{P}^d}(1)=-K_{X_d}$. It is easy to see that the degree of $X_d\subseteq\mathbb{P}^d$ is $H^2=d$. Further, $H\cdot(K_{X_d}+F)=-d+2$. We shall often denote the divisor $D=aL-b_1E_1-\cdots-b_{9-d}E_{9-d}$ by $(a;b_1,\cdots,b_{9-d})$.

\begin{theorem}\label{main'}
Let $X_d\subseteq\mathbb{P}^d$ be the anti--canonically embedded Del Pezzo surface of degree $3\leq d\leq 8$. Let $\mathcal{E}$ be a stable Ulrich vector bundle with respect to $H$ of rank $r$ where $H:=\mathcal{O}_{\mathbb{P}^d}(1)|_{X_d}=-K_{X_d}$. 
\begin{itemize}
    \item[(1)] If $d=3$ then $M_{\mathcal{E}}$ is simple and slope--semistable with respect to $H$. Moreover, the moduli spaces of slope--stable bundles  $\mathfrak{M}_H({\bf v}(\mathcal{E}))$ and $\mathfrak{M}_H({\bf v}(M_{\mathcal{E}}))$ are both smooth, irreducible and are birational to each other whenever the latter space is non--empty; in this case the general bundle in the latter moduli space is of the form $M_{\mathcal{F}}$ for some $\mathcal{F}\in \mathfrak{M}_{H}({\bf v}(\mathcal{E}))$.
    \item[(2)] If $d\geq 4$ then $\mathcal{S}_p(\mathcal{E})$ is simple and slope--semistable with respect to $H$ for all $p\geq 0$. Moreover, for any $p\geq 0$, the moduli space of slope--stable bundles $\mathfrak{M}_H({\bf v}(\mathcal{E}))$ and $\mathfrak{M}_H({\bf v}(\mathcal{S}_{p}(\mathcal{E})))$ are both smooth, irreducible and are birational to each other whenever the latter space is non--empty; in this case the general bundle in the latter moduli space is of the form $\mathcal{S}_p(\mathcal{F})$ for some $\mathcal{F}\in \mathfrak{M}_{H}({\bf v}(\mathcal{E}))$.
\end{itemize}
\end{theorem}

\noindent\textit{Proof.} First notice that for $d\geq 4$, the existence and slope--semistability of $\mathcal{S}_k(\mathcal{E})$ for all $k\geq 0$ follows from Theorem \ref{mains}. The simplicity for the iterated syzygy bundles in this case follows from Corollary \ref{extis}. When $d=3$, slope--semistability and simplicity for $M_{\mathcal{E}}$ follows from Proposition \ref{semistability} and Corollary \ref{subcan} respectively. Also, the smoothness and irreducibility of all the moduli spaces follow from Theorem \ref{wal}.

Now we prove the statement concerning birationality. Let $\mathbb{E}$ be a flat modular family of stable Ulrich bundles of rank and Chern classes same as $\mathcal{E}$ on $X\times S$ and denote by $p_1$ (resp. $p_2$) the projections onto $X$ (resp. $S$). Further let $\mathbb{H}:=p_1^*H$. Since $\mathbb{E}$ is relatively globally generated, the adjunction map $p_2^*{p_2}_*\mathbb{E}\to\mathbb{E}$ surjects and let $\mathbb{M}_{\mathbb{E}}$ be the kernel. Thus, we obtain the flat family $\mathbb{M}_{\mathbb{E}}$ on $X\times S/S$. Notice, for $s\in S$, $(\mathbb{M}_{\mathbb{E}})_s$ is $M_{\mathbb{E}_s}$ which is simple, slope--semistable and consequently $F$--prioritary. Thus, we obtain a morphism $S\to \mathfrak{Spl}_F({\bf v}(M_{\mathcal{E}}))$ that by Corollary \ref{subcan} descends to an \'etale dominant morphism $\mathfrak{M}_H({\bf v}(\mathcal{E}))\to\mathfrak{Spl}_F({\bf v}(M_{\mathcal{E}}))$. Now, when $\mathfrak{M}_H({\bf v}(M_{\mathcal{E}}))$ is non--empty, the restriction of the morphism to the pre--image of $\mathfrak{M}_H({\bf v}(M_{\mathcal{E}}))$ is injective by Corollary \ref{subcan} and thus $\mathfrak{M}_H({\bf v}(\mathcal{E}))$ and $\mathfrak{M}_H({\bf v}(M_{\mathcal{E}}))$ are birational. Notice that $\mathbb{M}_{\mathbb{E}}(\mathbb{H})$ is also a flat family on $X\times S/S$. An easy iteration of the above argument using Corollary \ref{extis} finishes the proof.\QEDB

\subsubsection{Moduli of stable bundles on cubic surfaces} Throughout this subsubsection, we work with anti--canonically embedded Del Pezzo surfaces $X_3\subseteq\mathbb{P}^3$. By \cite[Example 3.5]{CHGS}, the Ulrich line bundles on these surfaces correspond to twisted cubic curves. These curves are represented by 
\smallskip
\begin{center}
    $T_A=(1;0,0,0,0,0,0),\quad  T_B=(2;1,1,1,0,0,0),\quad T_C=(3;2,1,1,1,1,0)$ \\
    $T_D=(4;2,2,2,1,1,1),\quad T_E=(5;2,2,2,2,2,2)$
\end{center}
\smallskip
and their permutations; there are 72 of them. For future reference, we set $T_{B'}=(2;0,0,0,1,1,1)$ which is a permutation of the form listed for $T_B$. It is easy to see that $T^2=1$ for any twisted cubic curve.

\begin{lemma}\label{stablesum}
Let $\mathcal{E}$ be a stable Ulrich bundle with respect to $H$ on $X_3\subseteq\mathbb{P}^3$ of rank $r\geq 2$. Then there exist twisted cubic curves $T_1,\cdots ,T_r$ such that they satisfy the following two conditions simultaneously: 
 $$c_1(\mathcal{E})=\sum_{i=1}^r T_i,\quad \left(\sum\limits_{i=1}^{j-1} T_i\right)\cdot T_j\geq 2j-1\textrm{ for all }j=2,\cdots,r.$$

\end{lemma}
\noindent\textit{Proof.} For $r=2$, the statement follows from \cite[Theorem 4.3]{CHGS} and the fact that $T^2=1$ for any twisted cubic curve $T$. Assume the statements hold when $r\leq k-1$ and we prove this for $r=k\geq 3$. 
We have $c_1(\mathcal{E})T\geq 2k$ for any twisted cubic curve $T$ by \cite[Theorem 4.3]{CHGS}. Also notice that by the proof of \cite[Theorem 4.3]{CHGS}, there is a twisted cubic curve $T_k$ such that
$c_1(\mathcal{E})-T_k=c_1(\mathcal{F})$ where $\mathcal{F}$ is a stable Ulrich bundle of rank $k-1$. By induction hypothesis, there are twisted cubic curves $T_1,\cdots, T_{k-1}$ such that $c_1(\mathcal{F})=\sum_{i=1}^{k-1}T_i$ and $\left(\sum_{i=1}^{j-1} T_i\right)\cdot T_j\geq 2j-1$ for all $j=2,\cdots, k-1$. But by \cite[Theorem 4.3]{CHGS} again, we have
$$2k\leq \left(\sum_{i=1}^k T_i\right)\cdot T_k=\left(\sum_{i=1}^{k-1} T_i\right)\cdot T_k+T_k^2=\left(\sum_{i=1}^{k-1} T_i\right)\cdot T_k+1\implies \left(\sum_{i=1}^{k-1} T_i\right)\cdot T_k\geq 2k-1.$$
That completes the proof.\QEDB 

\begin{theorem}\label{cubic}
Let $\mathcal{E}$ be a stable Ulrich bundle with respect to $H$ on $X_3\subseteq\mathbb{P}^3$ of rank $r\geq 2$. Then the moduli spaces $\mathfrak{M}_H(r;c_1(\mathcal{E}),c_2(\mathcal{E}))$ and $\mathfrak{M}_H(2r;-c_1(\mathcal{E}), c_2(\mathcal{E})+r)$ are both non--empty, smooth, irreducible of dimension $(c_1(\mathcal{E}))^2-2r^2+1$. Moreover, they are birational to each other. 
\end{theorem}

\noindent\textit{Proof.} To start with, notice that $(c_1(\mathcal{E}))^2-c_2(\mathcal{E})=c_2(\mathcal{E})+r$ by \cite[Proposition 3.4]{CHGS} since $\mathcal{E}$ is Ulrich. Consequently, thanks to \eqref{cime}, Theorem \ref{main'} and Theorem \ref{wal}, it is enough to show that the moduli space $\mathfrak{M}_H(2r;-c_1(\mathcal{E}), c_2(\mathcal{E})+r)$ is non--empty.

By Lemma \ref{stablesum}, $c_1(\mathcal{E})=\sum_{i=1}^rT_i$ where the $T_i$'s are twisted cubic curves that satisfy the following 
\begin{equation}\label{decomp}
    \left(\sum\limits_{i=1}^{j-1} T_i\right)\cdot T_j\geq 2j-1\textrm{ for all }j=2,\cdots,r.
\end{equation}
It follows from \eqref{c1sum} and \eqref{cime} that 
\begin{equation}\label{eq0}
    c_1\left(\bigoplus_{i=1}^r M_{T_i}\right)=\sum_{i=1}^r c_1(M_{T_i})=-\sum_{i=1}^rT_i=-c_1(\mathcal{E}).
\end{equation}
We also compute the second Chern class of the direct sum using \eqref{c2sum} and \eqref{cime} to obtain
\begin{equation}\label{eq1}
    c_2\left(\bigoplus_{i=1}^r M_{T_i}\right)=\sum_{i=1}^r c_2(M_{T_i})+\sum_{i,j\in I_r,i<j}c_1(M_{T_i})c_1(M_{T_j})=\sum_{i=1}^rT_i^2+\sum_{i,j\in I_r,i<j} T_iT_j=r+\sum_{i,j\in I_r,i<j} T_iT_j.
\end{equation}
We use \cite[Proposition 3.4]{CHGS} again to obtain $c_2(\mathcal{E})+r=\frac{1}{2}((c_1(\mathcal{E}))^2+r)$ and consequently we get 
\begin{equation}\label{eq2}
    c_2(\mathcal{E})+r=\frac{1}{2}\left(\left(\sum_{i=1}^rT_i\right)^2+r\right)=\frac{1}{2}\left(\left(\sum_{i=1}^rT_i^2\right)+2\sum_{i,j\in I_r,i<j} T_iT_j+r\right)=r+\sum_{i,j\in I_r,i<j} T_iT_j.
\end{equation}
Observe that $h^0(T)=3$ for any twisted cubic curves since they are Ulrich, thanks to Proposition \ref{properties}, and consequently $\textrm{rank}(\bigoplus_{i=1}^r M_{T_i})=2r$. Thus, by \eqref{eq0} \eqref{eq1} and \eqref{eq2}, it is enough to prove that there exists a vector bundle $\mathcal{F}$ of rank $r$ with $c_i(\mathcal{F})=c_i(\bigoplus_{i=1}^rM_{T_i})$ for $i=1,2$ that is slope--stable with respect to $H$. We claim that for all $j=1,\cdots, r$, there exists vector bundle $\mathcal{F}_j$ of rank $2j$ that is slope--stable with respect to $H$ and satisfies $c_i(\mathcal{F}_j)=c_i(\bigoplus_{i=1}^jM_{T_i})$ for $i=1,2$. To prove this, we use induction; clearly the claim is true for $j=1$ since $M_{T}$ is slope--stable with respect to $H$ by Corollary \ref{coprimes}. Assume we have constructed such a vector bundle $\mathcal{F}_{j-1}$ for some $2\leq j\leq r-1$. Notice that $\mu_H(\mathcal{F}_{j-1})=\mu_H(M_{T_k})$ for all $1\leq k\leq r$. We aim to calculate the Euler characteristics of $\mathcal{F}_{j-1}^*\otimes M_{T_j}$ and $M_{T_j}^*\otimes \mathcal{F}_{j-1}$. We only do the computation for $\mathcal{F}_{j-1}^*\otimes M_{T_j}$, the other one is similar.

To calculate the first Chern class of $\mathcal{F}_{j-1}^*\otimes M_{T_j}$, we use \eqref{c1tensor} to obtain (using \eqref{cime})
\begin{equation}\label{eq3}
    c_1(\mathcal{F}_{j-1}^*\otimes M_{T_j})=-2c_1(\mathcal{F}_{j-1})+(2j-2)c_1(M_{T_j})=2\sum_{i=1}^{j-1}T_i-(2j-2)T_j.
\end{equation}
Observe that $c_1(\mathcal{F}_{j-1}^*\otimes M_{T_j})K_{X_3}=0$. Now we use \eqref{c2tensorv} and \eqref{cime} to obtain the second Chern class:
\begin{align}\label{eq4}
\begin{split}
    c_2(\mathcal{F}_{j-1}^*\otimes M_{T_j})  = & \binom{2j-2}{2}T_j^2+(2j-2)T_j^2+(2(2j-2)-1)\left(\sum_{i=1}^{j-1}T_i\right)\cdot(-T_j)+\left(\sum_{i=1}^{j-1}T_i\right)^2\\
    & +2\left(\left(\sum_{i=1}^{j-1}T_i^2\right)+\sum_{i,i'\in _{j-1},i<i'}T_iT_{i'}\right).
\end{split}
\end{align}
We use \eqref{eq3} and \eqref{eq4} to calculate the Euler characteristic of $\mathcal{F}_{j-1}^*\otimes M_{T_j}$ using \eqref{rr} to obtain 
\begin{align}
    \begin{split}
        \chi(\mathcal{F}_{j-1}^*\otimes M_{T_j})=& (4j-4)+\frac{4}{2}\left(\left(\sum_{i=1}^{j-1}T_i\right)-(j-1)T_j\right)^2-\frac{(2j-2)(2j-1)}{2}T_j^2+(4j-5)\left(\sum_{i=1}^{j-1}T_iT_j\right)\\
        & -\left(\left(\sum_{i=1}^{j-1}T_i^2\right)+2\left(\sum_{i,i'\in I_{j-1},i<i'}T_iT_{i'}\right)\right)-2\left(\left(\sum_{i=1}^{j-1}T_i^2\right)+\left(\sum_{i,i'\in I_{j-1},i<i'}T_iT_{i'}\right)\right).
    \end{split}
\end{align}
A straightforward computation yields $\chi(\mathcal{F}_{j-1}^*\otimes M_{T_j})=2(j-1)-\sum_{i=1}^{j-1}T_iT_j$. Since $\sum_{i=1}^{j-1}T_iT_j\geq 2j-1$ by \eqref{decomp}, we deduce the following inequality
\begin{equation}\label{eq6}
    \chi(\mathcal{F}_{j-1}^*\otimes M_{T_j})\leq -1.
\end{equation}

Now notice that $h^0(\mathcal{F}_{j-1}^*\otimes M_{T_j})=0$. Indeed, since $\mathcal{F}_{j-1}$ and $M_{T_j}$ are slope--stable with same slope, hence $h^0(\mathcal{F}_{j-1}^*\otimes M_{T_j})=0$ if $j\geq 3$. On the other hand, if $j=2$ and $h^0(\mathcal{F}_{j-1}^*\otimes M_{T_j})\neq 0$ then $\mathcal{F}_{1}\cong M_{T_2}$. Thus, $c_1(\mathcal{F}_1)c_1(M_{T_2})=1$. On the other hand by \eqref{decomp}, we have $c_1(\mathcal{F}_1)c_1(M_{T_2})=T_1\cdot T_2\geq 3$ which is a contradiction. Also, $h^2(\mathcal{F}_{j-1}^*\otimes M_{T_j})=0$. Indeed, $h^2(\mathcal{F}_{j-1}^*\otimes M_{T_j})=h^0(\mathcal{F}_{j-1}\otimes M_{T_j}^*\otimes K_{X_3})$ and the conclusion follows since $\mu_H(M_{T_j})>\mu_H(\mathcal{F}_{j-1}\otimes K_{X_3})$ since $M_{T_j}$ and $\mathcal{F}_{j-1}\otimes K_{X_3}$ are both slope--stable. Thus, by \eqref{eq6}, $h^1(\mathcal{F}_{j-1}^*\otimes M_{T_j})=-\chi(\mathcal{F}_{j-1}^*\otimes M_{T_j})\neq0$.
Similarly, we obtain $h^1(M_{T_j}^*\otimes \mathcal{F}_{j-1})\neq 0$.

Finally, since $\mu_H(\mathcal{F}_{j-1})=\mu_H(M_{T_j})$, the vector bundle $\mathcal{F}_{j-1}\oplus M_{T_j}$ is slope--semistable, which implies that is $F$--prioritary and hence unobstructed. Thus, $\mathcal{F}_{j-1}\oplus M_{T_j}$ deforms to slope--stable vector bundle $\mathcal{F}_j$ (see for example \cite[Corollary B.3]{Huy}) and by induction we are done.\QEDB 

\vspace{5pt}

We are now ready to provide the

\vspace{5pt}

\noindent\textit{Proof of Theorem \ref{main}.} Thanks to Theorem \ref{main'}, we only need to show that 
$\mathfrak{M}_H({\bf v}(M_{\mathcal{E}}))$ is non--empty when $d=3$. But this is a consequence of Theorem \ref{cubic} when $r\geq 2$, and Corollary \ref{coprimes} when $r=1$.\QEDB 

\vspace{5pt}

\noindent\textit{Proof of Corollary \ref{main3}.} It follows from the proof of Theorem \ref{main} (or rather Theorem \ref{main'}) that a general bundle in $\mathfrak{M}_{H}({\bf v}(\mathcal{S}_k(\mathcal{E})))$ is an iterated syzygy bundle if $\mathcal{S}_k(\mathcal{E})$ exists and the space is non--empty. The conclusion follows from Lemma \ref{elt}.\QEDB

\begin{corollary}\label{corcubic}
Let $X_3\subseteq\mathbb{P}^3$ be a cubic Del Pezzo surface. Then $\mathfrak{M}_{H}(4;c_1,c_2)$ are non-empty, smooth, irreducible and rational for the values described in the table below where $c_1'$ is arbitrary. We also provide the dimensions of these spaces.
\begin{table}[H]
    \centering
    \begin{tabular}{c|c|c}
        \hline
        $c_1$ & $c_2$ & $\dim \mathfrak{M}_{H}(4,c_1,c_2)$\\
        \hline\hline
        $-T_A-T_C+4c_1'$ & $6c_1'^2-3(T_A+T_C)c_1'+5$ & $1$\\
        \hline
        $-T_B-T_{B'}+4c_1'$ & $6c_1'^2-3(T_B+T_{B'})c_1'+6$ & $3$\\
        \hline
        $-T_A-T_{E}+4c_1'$ & $6c_1'^2-3(T_A+T_{E})c_1'+7$ & $5$\\
        \hline
    \end{tabular}
\end{table}
\end{corollary}
\noindent\textit{Proof.} Since the twist of a slope--stable vector bundle by a line bundle remains slope--stable, the statement is only about 
$\mathfrak{M}_H(4;(-4;-2,-1,-1,-1,-1,0),5)$, $\mathfrak{M}_H(4;(-4;-1,-1,-1,-1,-1,-1),6)$ and the moduli space $\mathfrak{M}_H(4;(-6;-2,-2,-2,-2,-2,-2),7)$. By Theorem \ref{cubic}, these spaces are irreducible, non--empty and birational to $\mathfrak{M}_H(2;(4;2,1,1,1,1,0),3)$, $\mathfrak{M}_H(2;(4;1,1,1,1,1,1),4)$ and  $\mathfrak{M}_H(2;(6;2,2,2,2,2,2),5)$ respectively, thanks to \cite[Example 3.6]{CHGS}. Twisting again by appropriate line bundles, the rationality of the latter spaces follow from \cite[Proposition 4.3.4 and Proposition 4.3.2]{CM2}.\QEDB

\subsubsection{Moduli of stable bundles on Del Pezzo surfaces of degree $d\geq 4$} Throughout this subsubsection, we work with anticanonically embedded Del Pezzo surfaces $X_d\subseteq\mathbb{P}^d$ for $4\leq d\leq 8$. Let $\mathcal{E}$ be an Ulrich bundle of rank $r$ with respect to $H$. We first compute the rank and the Chern classes of the bundles $\mathcal{S}_p(\mathcal{E})$. We set $N_{d,k}(\mathcal{E}):=\textrm{rank}(\mathcal{S}_k(\mathcal{E}))$. When the context is clear, we will abuse notation and simply write $N_{d,k}$ or $N_k$ for $N_{d,k}(\mathcal{E})$. Clearly $N_{-1}=r$ and $N_0=r(d-1)$.

\begin{lemma}\label{nk} The values of $N_k$ can be computed by the recurrence relation $N_k=(d-2)N_{k-1}-N_{k-2}$ for $k\geq 1$. In particular, we have the following formulae for all $k\geq -1$.
\begin{itemize}
    \item[(1)] When $d=4$, $N_k=(2k+3)r$.
    \item[(2)] When $d\geq 5$, set $\alpha_1:=(1/2)((d-2)+\sqrt{d(d-4)})$ and $\alpha_2:=(1/2)((d-2)-\sqrt{d(d-4)})$. Then
    $$N_k=\frac{r\left(\left(\alpha_2^{-(k+2)}+\alpha_2^{-(k+1)}\right)-\left(\alpha_1^{-(k+2)}+\alpha_1^{-(k+1)}\right)\right)}{\sqrt{d(d-4)}}.$$
\end{itemize}
In particular, $N_k\to\infty$ as $k\to\infty$.
\end{lemma}
\noindent\textit{Proof.} First we prove the recurrence relation. Notice that $\textrm{rank}(\mathcal{S}_1(\mathcal{E}))$ is the same as $\textrm{rank}(\mathcal{S}_1(\mathcal{E})|_C)$ where $C\in |H|$ is a smooth elliptic curve. By Proposition \ref{restricts}, $\mathcal{S}_1(\mathcal{E})|_C\cong \mathcal{S}_1(\mathcal{E}|_C)$ and thus by Lemma \ref{propertys} (2), $N_1=\textrm{rank}(M_{\mathcal{S}_0(\mathcal{E}|_C)})$ which in turn is equal to $h^0(\mathcal{S}_0(\mathcal{E}|_C))-\textrm{rank}(\mathcal{S}_0(\mathcal{E}|_C))=h^0(\mathcal{S}_0(\mathcal{E}|_C))-\textrm{rank}(\mathcal{S}_0(\mathcal{E}))$ where the last equality is obtained from another application of Proposition \ref{restricts}. We have the following short exact sequence
$$0\to \mathcal{S}_0(\mathcal{E}|_C)\to H^0(\mathcal{E}|_C)\otimes H|_C\to \mathcal{E}|_C(H|_C)\to 0.$$ Recall that (by Theorem \ref{mains}) $\mathcal{E}|_C$ satisfies $(\textrm{Kos}_{\infty})$ and $\mathcal{E}|_C$ is semistable of degree $rd$. Also, since $C$ is elliptic, for any semistable bundle $\mathcal{F}$ on $C$ with $\mu(\mathcal{F})\geq 2$, $h^0(\mathcal{F})=\textrm{deg}(\mathcal{F})$ by Riemann--Roch formula. Consequently we obtain $h^0(\mathcal{S}_0(\mathcal{E}|_C))=rd^2-2rd$. Thus $N_1=rd^2-2rd-r(d-1)=(d-2)N_0-N_{-1}$. 

Assume the recurrence relation is true for all $1\leq k\leq k_0$. As before, one computes
\begin{equation}\label{nk1}
    N_{k_0+1}=h^0(\mathcal{S}_{k_0}(\mathcal{E}|_C))-N_{k_0}.
\end{equation}
We have the following short exact sequence
$$0\to \mathcal{S}_{k_0}(\mathcal{E}|_C)\to H^0(\mathcal{S}_{k_0-1}(\mathcal{E}|_C))\otimes H|_C\to \mathcal{S}_{k_0-1}(\mathcal{E}|_C)(H|_C)\to 0.$$ 
Recall that $\mathcal{S}_p(\mathcal{E}|C)$ is semistable with slope $\geq 2$. Thus we obtain (applying Proposition \ref{restricts} if necessary)
\begin{equation}\label{nk2}
    N_{k_0}=h^0(\mathcal{S}_{k_0-1}(\mathcal{E}|_C))-N_{k_0-1},\quad h^0(\mathcal{S}_{k_0}(\mathcal{E}|_C))=d(h^0(\mathcal{S}_{k_0-1}(\mathcal{E}|_C))-h^0(\mathcal{S}_{k_0-1}(\mathcal{E}|_C)(H|_C)).
\end{equation}
We also have $h^0(\mathcal{S}_{k_0-1}(\mathcal{E}|_C)(H|_C))=h^0(\mathcal{S}_{k_0-1}(\mathcal{E}|_C))+dN_{k_0-1}$. Combining this with \eqref{nk1} and \eqref{nk2} we obtain 
$$N_{k_0+1}=d(N_{k_0}+N_{k_0-1})-(N_{k_0}+N_{k_0-1}+dN_{k_0-1})-N_{k_0}=(d-2)N_{k_0}-N_{k_0-1}.$$

It follows immediately from the recurrence relation that $N_k$'s are in arithmetic progression when $d=4$ and (1) follows. For (2), it is elementary to solve this recurrence. For simplicity, let us set $N_k=a_{k+1}$ for $k\geq -1$ and let $A(x)=\sum_{k=0}^{\infty}a_kx^k$. We have $$\sum_{k=2}^{\infty}(a_k-(d-2)a_{k-1}+a_{k-2})x^k=0\implies A(x)=\frac{r(x+1)}{x^2-(d-2)x+1}.$$ A straightforward computation yields that $$A(x)=\sum_{k=0}^{\infty}\frac{r\left(\left(\alpha_2^{-(k+1)}+\alpha_2^{-k}\right)-\left(\alpha_1^{-(k+1)}+\alpha_1^{-k}\right)\right)}{\sqrt{d(d-4)}}x^k$$ and the conclusion follows. Finally, by the recursive formula and induction, we obtain that $N_k-N_{k-1}>0$ for $k\geq 0$ and the last statement follows.\QEDB

\begin{corollary}\label{cornk}
The discriminant $\Delta(\mathcal{S}_k(\mathcal{E}))\to\infty$ as $k\to\infty$.
\end{corollary}

\noindent\textit{Proof.} Notice that $\Delta(\mathcal{S}_k(\mathcal{E}))=\dim\textrm{Ext}^1(\mathcal{S}_k(\mathcal{E}),\mathcal{S}_k(\mathcal{E}))+(N_k^2-1)$ which is equal to $\dim\textrm{Ext}^1(\mathcal{E},\mathcal{E})+(N_k^2-1)$ by Corollary \ref{extis}. The conclusion now follows from Lemma \ref{nk}.\QEDB

\vspace{5pt}

One can easily verify the following result by induction and we omit the proof.

\begin{lemma}\label{cink} As before, $\mathcal{E}$ is Ulrich on $X_d\subseteq\mathbb{P}^d$ with respect to $H$
 of rank $r$.
 \begin{itemize}
     \item[(1)] $c_1(\mathcal{S}_0(\mathcal{E})(-H))=-c_1(\mathcal{E})$ and for $k\geq 1$, we have the following formula
     $$c_1(\mathcal{S}_k(\mathcal{E})(-H))=(-1)^{k+1}c_1(\mathcal{E})+\sum_{i=0}^{k-1}(-1)^{k+i}N_{d,i}H.$$
     \item[(2)] $c_2(\mathcal{S}_0(\mathcal{E})(-H))=c_1(\mathcal{E})^2-c_2(\mathcal{E})$ and for $k\geq 1$, we have the following formula 
     \[
     \begin{array}{ll}
        c_2(\mathcal{S}_k(\mathcal{E})(-H))=  & \sum\limits_{i=0}^{k-1}(-1)^{k+i+1}\left(1-\binom{N_{d,i}}{2}\right)(c_1(\mathcal{S}_i(\mathcal{E})(-H)))^2+\\
        & \sum\limits_{i=0}^{k-1} (-1)^{k+i+1}(N_{d,i}+1)c_1(\mathcal{S}_i(\mathcal{E})(-H))H+\\
          & \sum\limits_{i=0}^{k-1}dN_{d,i}^2+(-1)^k(c_1(\mathcal{E})^2-c_2(\mathcal{E})).
     \end{array}\]
 \end{itemize}
\end{lemma}

\noindent\textit{Proof of Theorem \ref{main2}.} It follows directly from Theorem \ref{main}, Lemma \ref{nk}, Corollary \ref{cornk}, Lemma \ref{cink} combined with \cite[Proposition 5.1 (and Examples 6.4, 6.5, Theorem 6.7)]{Cas2}, and \cite[Propositions 4.3.2, 4.3.4 and 4.3.7]{CM2}.\QEDB

\vspace{5pt}

Finally we remark that in Theorem \ref{cubic} we showed the slope--stability of syzygy bundle of a stable Ulrich bundle $\mathcal{E}$ on cubic surfaces when the bundle is {\it general}. That motivates the following question which we believe has an affirmative answer, at least in the case of anti--canonical rational surfaces.
\begin{question}
Let $(X,H)$ be a polarized smooth projective variety with $H$ very ample. If $\mathcal{E}$ is a stable Ulrich bundle for $(X,H)$, is $M_{\mathcal{E}}$ slope--stable?
\end{question}

\section*{Declarations}

\subsection*{Conflict of interest} On behalf of all authors, the corresponding author states that there is
no conflict of interest.

\bibliographystyle{plain}

\end{document}